\documentclass[12pt]{amsart}

\usepackage{amscd,amssymb,amsthm,verbatim,
amsmath,amsfonts,mathrsfs,xspace,bbold}
\usepackage{color}

\newcommand{\Ad}{\operatorname{Ad}}
\newcommand{\ad}{\operatorname{ad}}

\newcommand{\dist}{\operatorname{dist}}
\renewcommand{\div}{\operatorname{div}}

\newcommand{\E}{\mathbb{E}}

\newcommand{\g}{\mathfrak{g}}

\renewcommand{\L}{{\mathscr L}}
\renewcommand{\H}{{\mathbb H}}
\renewcommand{\S}{{\mathscr S}}

\newcommand{\N}{{\mathbb N}}

\newcommand{\R}{{\mathbb R}}
\newcommand{\G}{{\mathbb G}}
\newcommand{\K}{{\mathbb K}}

\newcommand{\uno}{{\mathbb 1}}

\renewcommand{\span}{\operatorname{span}}

\newcommand{\vol}{\operatorname{vol}}

\newcommand{\eps}{\epsilon}

\usepackage[all]{xy}

\def\XXint#1#2#3{{\setbox0=\hbox{$#1{#2#3}{\int}$}
\vcenter{\hbox{$#2#3$}}\kern-.5\wd0}}

\numberwithin{equation}{section}
\theoremstyle{plain}

\newtheorem{lemma}[equation]{Lemma}
\newtheorem{theorem}[equation]{Theorem}
\newtheorem{proposition}[equation]{Proposition}

\newtheorem{definition}[equation]{Definition}

\theoremstyle{remark}
\newtheorem{remark}[equation]{Remark}
\begin{document}
\title[Tangent hyperplane in Carnot groups]
{Rectifiability of sets of finite perimeter in \\ Carnot groups:
existence of a tangent hyperplane}
\author{Luigi Ambrosio}
\address{Scuola Normale Superiore, Pisa}
\email{l.ambrosio@sns.it}
\author{Bruce Kleiner}
\address{Yale University, USA}
\email{bruce.kleiner@yale.edu}
\author{Enrico Le Donne}
\address{Yale University, USA}
\email{enrico.ledonne@yale.edu}
\thanks{The second author was partially supported by NSF grant 
DMS-0701515}
\date{\today}

\begin{abstract}
We consider sets of locally finite perimeter in Carnot
groups.  We show that if $E$ is a set of locally finite perimeter
in a Carnot group $G$ then, for almost every $x\in G$ with respect
to the perimeter measure of $E$, some tangent of $E$ at $x$
is a vertical halfspace. This is a partial extension of 
a theorem of Franchi-Serapioni-Serra Cassano 
in step 2 Carnot groups: they show in \cite{fsc2,fsc} that, for almost every $x$,
$E$ has a unique tangent at $x$, and this tangent is a vertical halfspace.
\end{abstract}
\maketitle

\section{Introduction}

The differentiability properties of functions and the rectifiability
properties of sets are classical themes of Real Analysis and
Geometric Measure Theory, with many mutual connections. In the
context of stratified Carnot groups, the first problem has been
solved, within the category of Lipschitz maps, in a deep work of
Pansu \cite{pansu}; here we are interested in the second problem, in
the class of sets $E$ of locally finite perimeter: if we denote by
$X_1,\ldots,X_m$ an orthonormal basis of the horizontal layer of the
Lie algebra $\g$ of left-invariant vector fields of the Carnot group
$\G$, this class of sets is defined by the property that the
distributional derivatives $X_1\uno_E,\ldots,X_m\uno_E$ are
representable by Radon measures in $\G$. This notion, which extends
the classical one developed and deeply studied by De Giorgi in
\cite{dg1} and \cite{dg2} (see also \cite{Ambook}), is compatible
with the Carnot-Carath\'eodory (subriemannian) distance $d$ induced
by $X_1,\ldots,X_m$; in this context the total variation $|D\uno_E|$
of the $\R^m$-valued measure $(X_1\uno_E,\ldots,X_m\uno_E)$ plays
the role of surface measure associated to $d$. Our interest in this
topic was also motivated by the recent papers \cite{cheekle1},
\cite{cheekle2}, where sets of finite perimeter in Carnot groups
(and in particular in the Heisenberg groups) are used to study a new
notion of differentiability for maps with values in $L^1$, with the
aim of finding examples of spaces which cannot be bi-Lipschitz
embedded into $L^1$.

The first basic properties of the class of sets of finite perimeter
(and of $BV$ functions as well),
such as compactness, global and local isoperimetric inequalities,
have been proved in \cite{garnieu}; then, in a series of papers
\cite{fsc2,fsc}, Franchi, Serapioni and Serra Cassano made a more
precise analysis of this class of sets, first in the Heisenberg
groups $\H^n$ and then in all step 2 Carnot groups (using also some
measure-theoretic properties proved, in a more general context,
in \cite{Amb01}, see also Theorem~\ref{tambrosio}). As in the work
of De Giorgi, the crucial problem is the analysis of tangent sets to
$E$ at a point $\bar x$, i.e. all limits
$$
\lim_{i\to\infty}\delta_{1/r_i}(\bar x^{-1}E)
$$
where $(r_i)\downarrow 0$ and convergence occurs locally in measure
(here $\delta_r:\G\to\G$ denote the intrinsic dilations of the
group). In \cite{fsc} it is proved that for $|D\uno_E|$-a.e. $\bar
x$ there exists a unit vector $\nu_E(\bar x)\in{\bf S}^{m-1}$, that
we shall call horizontal normal, such that
\begin{equation}\label{invia}
\sum_{i=1}^m\nu_{E,i}(\bar x)X_i\uno_E\geq 0\qquad\text{and}\qquad
\sum_{i=1}^m\xi_i(\bar x)X_i\uno_E=0\quad\forall\xi\perp\nu_E(\bar
x).
\end{equation}
We shall call these sets
with constant horizontal normal (identified, in the coordinates
relative to the basis $X_1,\ldots,X_m$, by the vector $\nu_E(\bar
x)$): the question is whether \eqref{invia} implies additional
information on the derivative of $E$ along vector fields $Y$ that
do not belong to the horizontal layer: even though $m<n={\rm
dim}(\g)$, this can be expected, having in mind that the Lie algebra
generated by $X_1,\ldots,X_m$ is the whole of $\g$. The main result
of \cite{fsc} is the proof that, in all step 2 groups, \eqref{invia}
implies $[X_i,X_j]\uno_E=0$ for all $i,\,j=1,\ldots,m$. As a
consequence, up to a left translation $E$ is really, when seen in
exponential coordinates, an halfspace:
$$
\left\{x\in\R^n:\ \sum_{i=1}^m\nu_{E,i}(\bar x)x_i\geq 0\right\}.
$$
We shall call it vertical halfspace, keeping in mind that there is
no dependence on the coordinates $x_{m+1},\ldots,x_n$. This fact
leads to a complete classification of the tangent sets and has
relevant consequences, as in the classical theory, on the
representation of $|D\uno_E|$ in terms of the spherical Hausdorff
measure and on the rectifiability, in a suitable intrinsic sense, of
the measure-theoretic boundary of $E$, see \cite{fsc} for more
precise informations.

On the other hand, still in \cite{fsc}, it is proved that for
general Carnot groups the conditions \eqref{invia} do not
characterize vertical halfspaces: an explicit example is provided in
a step 3 group of Engel type (see also Section~\ref{sengel}). Basically,
because of this obstruction, the results of \cite{fsc} are
limited to step 2 groups.

The classification and even the regularity properties of sets $E$
with a constant horizontal normal is a challenging and, so far,
completely open question. However, recently we found a way to bypass
this difficulty and, in this paper, we show the following result:

\begin{theorem}
\label{main} Suppose $E\subseteq\G$ has locally finite perimeter.
Then, for $|D\uno_E|$-a.e. $\bar x\in\G$ a vertical halfspace $H$
belongs to the tangents to $E$ at $\bar x$.
\end{theorem}

Of course Theorem~\ref{main} does not provide yet a complete
solution of the rectifiability problem: indeed, even though the
direction $\nu_E(\bar x)$ of the halfspace $H$ depends on $\bar x$
only, we know that $\bar x^{-1}E$ is close on an infinitesimal
sequence of scales to $H$, but we are not able to show that this
happens on all sufficiently small scales. What is still missing is
some monotonicity/stability argument that singles out halfspaces as
the only possible tangents, wherever they are tangent (see also the
discussion in Remark~\ref{rloca}). In a similar context, namely the
rectifiability of measures having a spherical density, this is
precisely the phenomenon discovered by Preiss in \cite{Preiss}: we
took some ideas from this paper, adapting them to the setting of
Carnot groups, to obtain our result. For these reasons, the complete
solution of the rectifiability problem seems to be related to the
following question (we denote by $\vol_\G$ the Haar measure of the group 
and by $e$
the identity of the group): let $E\subset\G$ be a set with a constant
horizontal normal $\nu\in {\bf S}^{m-1}$ and let $H$ be a vertical
halfspace with the same horizontal normal; if 
$$
\liminf_{R\to+\infty}\frac{\vol_\G\bigl((E\Delta H)\cap B_{
R}(e)\bigr)}{\vol_\G\bigl(B_{R}(e)\bigr)}=0,
$$
is it true that $E$ is a vertical halfspace?
However, as pointed out to us by Vittone, 
the answer to this question is negative,
see \eqref{toobad}, so that new ideas seem to be needed to
prove the uniqueness, at $|D\uno_E|$-a.e. point, of the tangent set.

In order to illustrate the main ideas behind the proof of our result,
let us call regular directions of $E$ the vector fields $Z$ in the
Lie algebra $\g$ such that $Z\uno_E$ is representable by a Radon
measure, and invariant directions those for which the measure is
$0$. Our strategy of proof rests mainly on the following
observations: the first one (Proposition~\ref{pollo1}) is that
the adjoint operator $\Ad_{\exp(Y)}:\g\to\g$ maps regular directions into regular
directions whenever $Y$ is an invariant direction. If
$$
X:=\sum_{i=1}^m\nu_{E,i}(\bar x)X_i\in\g,
$$
we look at the vector space spanned by $\Ad_{\exp(Y)}(X)$, as $Y$
varies among the invariant directions, and use this fact to show
that any set with constant horizontal normal must have a regular
direction $Z$ not belonging to the vector space spanned by the
invariant directions and $X$ (which contains at least the horizontal
layer). This is proved in Proposition~\ref{pcrucial} in purely
geometric terms in general Lie gropus, and
Proposition~\ref{pcrucialbis} provides a more explicit expression of
the new regular directions generated, in Carnot groups, with this
procedure.

Then, the second main observation is that if a regular direction $Z$
for a set $F$ has no component in the horizontal layer, then the
tangents to $F$ at $\bar x$ are invariant along a new direction
depending on $Z$ for most points $\bar x$; this follows
(Lemma~\ref{provafis}) by a simple scaling argument, taking into
account that the Lie algebra dilations $\delta_r$ shrink more, as
$r\downarrow 0$, in the non-horizontal directions. Therefore, at
many points, a tangent to a set with constant horizontal normal has
a new invariant direction. Having gained this new direction, this
procedure can be restarted: the adjoint can be used to generate a
new regular direction, then a tangent will have a new invariant
direction, and so on.

In this way we show in Theorem~\ref{main1} that, if we iterate the
tangent operator sufficiently many times (the number depending on
the Lie algebra stratification only) we do get a vertical halfspace.
This means that we consider a tangent set $E^1$ to $E$ at $\bar x$,
then a tangent $E^2$ to $E^1$ at a suitable point $\bar x_1$ in the
support of $|D\uno_{E^1}|$, and so on. At this stage we borrow some
ideas from \cite{Preiss} to conclude that, at $|D\uno_E|$-a.e. point
$\bar x$, iterated tangents are tangent to the initial set: this is
accomplished in Section~\ref{sitera} and leads to the proof of
Theorem~\ref{main}.

\medskip
{\bf Acknowledgements.} We thank V. Magnani and A. Martini for some
useful comments on a preliminary version of this paper.

\section{Main notions}

\subsection{Vector fields, divergence, $X$-derivative}
Throughout this section, we will denote by $M$ a smooth
differentiable manifold with topological dimension $n$, endowed with
a $n$-differential volume form $\vol_M$ (eventually $M$ will be a
Lie group $\G$, and $\vol_M$ the right Haar measure).

For $x\in M$, the fiber $T_xM$ of the tangent bundle $TM$
is a derivation of germs of $C^\infty$ functions at $x$
(i.e., an $\R$-linear application
 from $C^\infty(x)\to\R$ that satisfies the Leibnitz rule).
If $F:M\to N$ is smooth and $x\in M$,
we shall denote by $dF_x:T_xM\to T_{F(x)}N$
its differential, defined as follows: the pull back operator
$u\mapsto F_x^*(u):=u\circ F$ maps $C^\infty\left(F(x)\right)$
into $C^\infty(x)$; thus, for $v\in T_xM$ we have that
$$
dF_x(v)(u):=v(u\circ F)(x),
\qquad u\in C^\infty(F(x))
$$
defines an element of $T_{F(x)}N$.

We denote by $\Gamma(TM)$ the linear space of smooth vector fields,
i.e. smooth sections of the tangent bundle $TM$; we will typically
use the notation $X,\,Y,\,Z$ to denote them. We use the notation
$[X,Y]f:=X(Yf)-Y(Xf)$ for the Lie bracket, that induces on
$\Gamma(TM)$ an infinite-dimensional Lie algebra structure.

If $F:M\to N$ is smooth and invertible and $X\in\Gamma(TM)$,
the push forward vector field $F_*X\in \Gamma(TN)$ is defined by
the identity
$(F_*X)_{F(x)}=dF_x(X_x)$. Equivalently,
\begin{equation} (F_*X)u:=[X(u\circ F)]\circ F^{-1}\qquad\forall
u\in C^\infty(M).\end{equation}
The push-forward commutes with the Lie bracket, namely
\begin{equation}\label{bir}
[F_*X,F_*Y]=F_*[X,Y]\qquad\forall X,\,Y\in\Gamma(TM).
\end{equation}

If $F:M\to N$ is smooth and $\sigma$ is a smooth curve on $M$, then
\begin{equation}
dF_{\sigma(t)}(\sigma'(t))=(F\circ\sigma)'(t),\end{equation} where
$\sigma'(t)\in T_{\sigma(t)}M$ and $(F\circ\sigma)'(t)\in
T_{F(\sigma(t))}N$ are the tangent vector fields along the two
curves, in $M$ and $N$. If $u\in C^\infty(M)$, identifying
$T_{u(p)}\R$ with $\R$ itself, given $X\in \Gamma(TM)$, we have
$$du_p(X)=X_p(u).$$

Now we use the volume form to define the divergence as follows:
\begin{equation}\label{def_divergence}
\int_M X u\,d\vol_M=-\int_M u \div X\,d\vol_M\qquad\forall u\in
C^\infty_c(M).
\end{equation}
When $(M,g)$ is a Riemannian manifold and $\vol_M$ is the volume
form induced by $g$, then an explicit expression of this differential
operator can be obtained in terms of the components of $X$, and
\eqref{def_divergence} corresponds to the divergence theorem
on manifolds. We
won't need either a Riemannian structure or an explicit expression
of $\div X$ in the sequel, and for this reason we have chosen a
definition based on \eqref{def_divergence}: this emphasizes the
dependence of $\div X$ on $\vol_M$ only. By applying this identity
to a divergence-free vector field $X$, we obtain
\begin{equation}\label{int_parts}
\int_M  u X v\,d\vol_M=-\int_M v Xu\,d\vol_M \qquad\forall u,\,v\in
C^\infty_c(M).
\end{equation}

This motivates the following classical definition.

\begin{definition}[X-distributional derivative]
Let $u\in L^1_{\rm loc}(M)$ and let $X\in\Gamma(TM)$ be divergence-free.
We denote by $Xu$ the distribution
$$
\langle Xu,v\rangle:=-\int_M uXv\,d\vol_M,\qquad v\in C^\infty_c(M).
$$
If $f\in L^1_{\rm loc}(M)$, we write $Xu=f$ if $\langle
Xu,v\rangle=\int_M vf\,d\vol_M$ for all $v\in C^\infty_c(M)$.
Analogously, if $\mu$ is a Radon measure in $M$, we write $Xu=\mu$
if $\langle Xu,v\rangle=\int_M v\,d\mu$ for all $v\in
C^\infty_c(M)$.
\end{definition}

According to \eqref{int_parts} (still valid when $u\in C^1(M)$), the
distributional definition of $Xu$ is equivalent to the classical one
whenever $u\in C^1(M)$.

In Euclidean spaces, the $X$-derivative of characteristic functions
of nice domains can be easily computed (and of course the result
could be extended to manifolds, but we won't need this extension).

\subsection{$X$-derivative of nice functions and domains}\label{rgauss}
If $u$ is a $C^1$ function in $\R^n$, then $Xu$ can be calculated as
the scalar product between $X$ and the gradient of $u$:
\begin{equation}\label{prodgrad}Xu=\langle X,\nabla u\rangle.\end{equation}
Assume that $E\subset\R^n$ is locally the sub-level set
of the $C^1$ function $f$ and that $X\in\Gamma(T\R^n)$ is
divergence-free. Then, for any $v\in C^\infty_c(\R^n)$ we can apply
the Gauss--Green formula to the vector field $v X$, whose divergence
is $Xv$, to obtain
$$
\int_E Xv\,dx=\int_{\partial E}\langle v X,\nu^{eu}_E\rangle\, d{\mathscr
H}^{n-1},
$$
where $\nu^{eu}_E$ is the unit (Euclidean) outer normal to $E$. This proves that
$$
X\uno_E=-\langle X,\nu^{eu}_E\rangle {\mathscr H}^{n-1}\llcorner_{\partial E}.
$$
However, we have an explicit formula for the unit (Euclidean)
outer normal to $E$, it is
$\nu^{eu}_E(x)=\nabla f(x)/|\nabla f(x)|$, so, by (\ref{prodgrad}),
\begin{eqnarray*}\langle X,\nu^{eu}_E\rangle &=&\langle X,\frac{\nabla f}{|\nabla f|}\rangle \\
&=& \frac{\langle X,\nabla f\rangle }{|\nabla f|} = \frac{Xf}{|\nabla f| }.
\end{eqnarray*}
Thus
\begin{equation}\label{Euclid}X\uno_E=- \frac{Xf}{|\nabla f| }  {\mathscr H}^{n-1}\llcorner_{\partial E}.
\end{equation}
\subsection{Flow of a vector field}

Given $X\in\Gamma(TM)$ we can consider the associated flow, i.e.,
the solution $\Phi_X:M\times \R\to M$ of the following ODE
\begin{equation}\label{flow}
 \left\{ \begin{array}{ccl}
\displaystyle{\frac{d}{dt}}\Phi_X(p,t)&=&X_{\Phi_X(p,t)} \\ \\
\Phi_X(p,0)&=&p.
\end{array} \right.
\end{equation}
Notice that the smoothness of $X$ ensures uniqueness, and therefore
the semigroup property
\begin{equation}\label{semigroup}
\Phi_X(x,t+s)=\Phi_X(\Phi_X(x,t),s)\qquad\forall t,\,s\in
\R,\,\,\forall x\in M
\end{equation}
 but not global existence; it will be guaranteed,
however, in all cases considered in this paper. We obviously have
\begin{equation}\label{esns1}
\frac{d}{d t}(u\circ\Phi_X)(p,t)=(Xu)(\Phi_X(p,t)) \qquad\forall
u\in C^1(M).
\end{equation}
An obvious consequence of this identity is that, for a $C^1$
function $u$, $Xu=0$ implies that $u$ is constant along the flow,
i.e. $u\circ\Phi_X(\cdot,t)=u$ for all $t\in\R$.  A similar
statement holds even for distributional derivatives along vector
fields: for simplicity let us state and prove this result for
divergence-free vector fields only.

\begin{theorem}\label{tdp} Let $u\in L^1_{\rm loc}(M)$
be satisfying $Xu=0$ in the sense of distributions. Then, for all
$t\in\R$, $u=u\circ\Phi_X(\cdot,t)$ $\vol_M$-a.e. in $M$.
\end{theorem}
\begin{proof} Let $g\in C^1_c(M)$; we need to show that the map
$t\mapsto\int_M g u\circ\Phi_X(\cdot,t)\,d\vol_M$ is independent of
$t$. Indeed, the semigroup property \eqref{semigroup}, and the fact
that $X$ is divergence-free yield
\begin{eqnarray*}
&&\int_M g u\circ\Phi_X(\cdot,t+s)\,d\vol_M- \int_M g
u\circ\Phi_X(\cdot,t)\,d\vol_M \\&=& \int_M u
g\circ\Phi_X(\cdot,-t-s)\,d\vol_M- \int_M u
g\circ\Phi_X(\cdot,-t)\,d\vol_M \\&=&\int_M u
g\circ\Phi_X(\Phi_X(\cdot,-s),-t)\,d\vol_M-\int_M u
g\circ\Phi_X(\cdot,-t)\,d\vol_M\\&=& -s\int_M u
X(g\circ\Phi_X(\cdot,-t))\,d\vol_M+o(s)=o(s).
\end{eqnarray*}
\end{proof}

\begin{remark}\label{renrico}{\rm
We notice also that the flow is $\vol_M$-measure preserving (i.e.
$\vol_M(\Phi_X(\cdot,t)^{-1}(A))=\vol_M(A)$ for all Borel sets
$A\subseteq M$ and $t\in\R$) if and only if $\div X$ is equal to 0.
Indeed, if $f\in C^1_c(M)$, the measure preserving property gives
that $\int_M f(\Phi_X(x,t))\,d\vol_M(x)$ is independent of $t$. A
time differentiation and \eqref{esns1} then give
$$
0=\int_M\frac{d}{dt} f(\Phi_X(x,t))\,d\vol_M(x)= \int_M
Xf(\Phi_X(x,t))\,d\vol_M(x)=\int_M Xf(y)\,d\vol_M(y).
$$
Therefore $\int_M f\div X\,d\vol_M=0$ for all $f\in C^1_c(M)$, and
$X$ is divergence-free. The proof of the converse implication is
similar, and analogous to the one of Theorem~\ref{tdp}.}\end{remark}

\subsection{Lie groups}
Let $\G$ be a Lie group, i.e. a differentiable $n$-dimensional
manifold with a smooth group operation. We shall denote by $e$ the
identity of the group, by $R_g(h):=hg$ the right translation, and by
$L_g(h):=gh$ the left translation. We shall also denote by $\vol_\G$
the volume form and, at the same time, the right-invariant Haar
measure.

Forced to make a choice, we follow the majority of the literature
focusing on the \emph{left} invariant vector fields. i.e. the vector
fields $X\in \Gamma(T\G)$ such that $(L_g)_*X=X$, so that
$(dL_g)_xX=X_{L_g(x)}$ for all $x\in \G$. In differential terms, we
have
$$
X (f\circ L_g)(x)=Xf(L_g(x))\qquad\forall x,\,g\in \G.
$$
Thanks to \eqref{bir} with $F=L_g$, the class of left invariant
vector fields is easily seen to be closed under the Lie bracket, and
we shall denote by $\g\subseteq\Gamma(T\G)$ the \emph{Lie algebra}
of left invariant vector fields. We will typically use the notations
$U,\,V,\,W$ to denote subspaces of $\g$.

Note that, after fixing a vector $v\in T_e\G$, we can construct a
left invariant vector field $X$ defining $X_g:=(L_g)_*v$ for any
$g\in \G$. This construction is an isomorphism between the set $\g$
of all left invariant vector fields and $T_e\G$, and proves that
$\g$ is a $n$-dimensional subspace of $\Gamma(T\G)$.

Let $X\in\g$ and let us denote, as usual in the theory, by
$\exp(tX)$ the flow of $X$ at time $t$ starting from $e$ (that is,
$\exp(tX):=\Phi_X(e,t)=\Phi_{tX}(e,1)$); then, the curve $g\exp(tX)$
is the flow starting at $g$: indeed, since $X$ is left invariant,
setting for simplicity $\gamma(t):=\exp(tX)$ and
$\gamma_g(t):=g\gamma(t)$, we have
$$\frac{d}{dt}\gamma_g(t)=\frac{d}{dt}( L_g(\gamma(t)))=(dL_g)_{\gamma(t)}
\frac{d}{dt}\gamma(t)=(dL_g)_{\gamma(t)} X=X_{\gamma_g(t)}.$$

This implies that $\Phi_X(\cdot,t)=R_{\exp(t X)}$ and so the flow
preserves the right Haar measure, and the left translation preserves
the flow lines. By Remark~\ref{renrico} it follows that all $X\in\g$
are divergence-free, and Theorem~\ref{tdp} gives
\begin{equation}\label{righto}
f\circ R_{\exp(tX)}=f\quad\forall
t\in\R\qquad\Longleftrightarrow\qquad X f=0
\end{equation}
whenever $f\in L^1_{\rm loc}(\G)$.

Before stating the next proposition, we recall the definition of the
adjoint. For $k\in\G$, the conjugation map
\begin{eqnarray}
C_k:&\G\to& \G\nonumber\\
&g\mapsto& C_k(g):=kgk^{-1}
\end{eqnarray}
is the composition of $L_k$ with $R_{k^{-1}}$. The \emph{adjoint}
operator $k\mapsto\Ad_k$ maps $\G$ in $GL(\g)$ as follows:
\begin{eqnarray}
\text{$\Ad_k(X):= (C_k)_*X$, so that $\Ad_k(X)f(x)=X(f\circ
C_k)(C_k^{-1}(x))$.}
\end{eqnarray}
The definition is well posed because $\Ad_k(X)$ is left invariant
whenever $X$ is left invariant: for all $g\in \G$ we have indeed
$$
\Ad_k(X)(f\circ L_g)(x)=X(f\circ L_g\circ C_k)(k^{-1}xk)= X(f\circ
R_{k^{-1}}\circ L_{gk})(k^{-1}xk)=X(f\circ R_{k^{-1}})(gxk).
$$
On the other hand
$$
\Ad_k(X) f(L_g(x))=X(f\circ C_k)(k^{-1}gxk)=X (f\circ
R_{k^{-1}})(gxk).
$$

\begin{proposition}\label{pcrucial}
Assume that $\G$ is a connected, simply connected nilpotent Lie
group. Let $\g'$ be a Lie subalgebra of $\g$ satisfying ${\rm
dim}(\g')+2\leq{\rm dim}(\g)$, and assume that $W:=\g'\oplus\{\R
X\}$ generates the whole Lie algebra $\g$ for some $X\notin\g'$.
Then, there exists $k\in\exp(\g')$ such that $\Ad_k(X)\notin W$.
\end{proposition}
\begin{proof} Note that $\g'$ is a finite-dimensional sub-algebra
and that $\exp$ is, under the simple connectedness assumption, a
homeomorphism, hence $\K:=\exp(\g')$ is a closed (proper) Lie
subgroup of $\G$. Therefore, we can consider the quotient manifold
$\G/\K$, in fact the homogeneous space of right cosets: it consists
of the equivalence classes of $\G$ induced by the relation
$$
x\sim y\qquad\Longleftrightarrow\qquad y^{-1}x\in\K.
$$
We shall denote by $\pi:\G\to \G/\K$ the canonical projection. The natural topology of $\G/\K$ is determined
by the requirement that $\pi$ should be continuous and open. Let $\mathfrak{m}$ denote some vector space of $\g$
such that $\g=\g'\oplus\mathfrak{m}$. The sub-manifold $\exp(\mathfrak{m})$
is referred as a local cross section for $\K$ at the origin, and it can be used to
give a differentiable structure to $\G/\K$. In fact, let $Z_1,\ldots,Z_r$ be a basis
of $\mathfrak{m}$, then the mapping $$(x_1, \ldots,x_r)\mapsto \pi(g\exp(x_1Z_1+\ldots+x_rZ_r))$$
is a homeomorphism of an open set of $\R^r$ onto a neighborhood of $g\K$ in $\G/\K$.
Then it is easy (see \cite{Helgason} for details) to see that with these charts, $\G/\K$ is an analytic manifold.
In particular, $\pi$ restrict to $\exp(\mathfrak{m})$ is a local diffeomorphism into $\G/\K$
and $d\pi(X)\neq 0$ since the projection of $X$ on $\mathfrak{m}$ is non zero.

Notice that, by our assumption on the dimension of $\g'$, the
topological dimension of $\G/\K$ is at least 2. Now, if the
statement were false, taking into account that
$\Ad_k(\g')\subseteq\g'$, we would have $\Ad_k(W)\subseteq W$ for
all $k\in\K$. By the definition of adjoint operator as composition
of the differentials of right and left translations, the above would
be equivalent to say that
$$(R_k)_*((L_{k^{-1}})_*(Y))\in W\qquad\forall Y\in W,\,\,k\in\K.$$
Since the vector fields in $W$ are left invariant (i.e. $(L_g)_*Y=Y$
for all $Y\in W$), this condition would say that $W$ is $\K$-right
invariant, and we can write this condition in the form
$d(R_k)_x(W_x)\subset W_{xk}$ for all $x\in \G$ and $k\in \K$.

Now, let us consider the subspaces $d\pi_x(W_x)$ of $T_{\pi(x)}\G/\K$:
they are all 1-dimensional, thanks to the fact that ${\rm
dim}(W)=1+{\rm dim}(\g')$, and they depend only on $\pi(x)$: indeed,
$\K$-right invariance and the identity $\pi\circ R_k=\pi$ give
$$
d\pi_x(Y_x)=d\pi_{xk}(d(R_k)_x(Y_x))\in d\pi_{xk}(W_{xk})
$$
for all $Y\in W$ and $k\in\K$. Therefore we can define a (smooth)
1-dimensional distribution $W/\K$ in $\G/\K$ by
$(W/\K)_y:=d\pi_x(W_x)$, where $x$ is any element of $\pi^{-1}(y)$.
In particular $W/\K$ would be tangent to a $1$-dimensional foliation
${\mathcal F}$ of $\G/\K$ that has at least codimension $1$, since
$\G/\K$ has at least dimension $2$. Letting ${\mathcal F}'$ be the
foliation of $\G$ whose leaves are the inverse images via $\pi$ of
leaves of ${\mathcal F}$, we find that still ${\mathcal F}'$ has
codimension at least 1, and $W$ is tangent to the leaves of
${\mathcal F}'$. But this contradicts the fact that $W$ generates
$\g$: in fact, the only sub-manifold to which $W$ could be tangent
is all the manifold $\G$.
\end{proof}

In the following proposition we provide a characterization of the
vector space spanned by $\Ad_{\exp(Y)}(X)$, where $Y$ varies in a
Lie subalgebra of $\g$.  This improved version of
Proposition~\ref{pcrucialbis} was pointed out to us by V. Magnani.

\begin{proposition}\label{pcrucialbis} Let $\g$ be a nilpotent Lie algebra,
let $\g'\subset\g$ be a Lie algebra and let $X\in\g$. Then
$$
{\rm span}\left(\{\Ad_{\exp(Y)}(X):\ Y\in\g'\}\right)= [\g',X]+
[\g',[\g',X]]+[\g',[\g',[\g',X]]]+\cdots.
$$
\end{proposition}
\begin{proof}
Let us denote by $S$ the space $ {\rm
span}\left(\{\Ad_{\exp(Y)}(X):\ Y\in\g'\}\right)$. Obviously $S$
contains $X$ and all vector fields $\Ad_{\exp(rY)}(X)$ for $r\geq 0$
and $Y\in\g'$. Now, denoting by $L(\g)$ the linear maps from $\g$ to
$\g$, let us recall the formula (see \cite{knapp}, page 54)
$\Ad_{\exp(Y)}=e^{\ad_Y}$, where $\ad_\cdot:\g\to {\rm End}(\g)$ is
the operator $\ad_Y(X)=[Y,X]$ and the exponential $e^A$ is defined
for any $A\in L(\g)$, by $e^A:=\sum\limits_{i=0}^\infty A^i/i!\in
L(\g)$. Therefore
\begin{equation}\label{eqknapp}
\Ad_{\exp(Y)}X=X+[Y,X]+\dfrac{1}{2}[Y,[Y,X]]+\cdots
\end{equation}
Let $\nu$ be the dimension of $\g'$ and let $(Y_1,\ldots,Y_\nu)$ be
a basis of $\g'$. Taking into account the identity (\ref{eqknapp}),
for all $Y=\sum_1^\nu r_j Y_j\in\g'$, we define
\begin{eqnarray*}
\Phi(r_1,\ldots,r_\nu)&:=&\Ad_{\exp(\sum_1^\nu r_j Y_j)}X-X\\
&=&\sum_{k=1}^{s-1}\frac{1}{k!}
\big(\sum_{j=1}^\nu r_j\, \ad Y_j\big)^kX\\
&=&\sum_{k=1}^{s-1}\frac{1}{k!} \sum_{j_1,\ldots,j_k=1}^\nu
r_{j_1}\cdots r_{j_k} \big(\ad Y_{j_1}\cdots\ad Y_{j_k}\big)X\in
S\,.
\end{eqnarray*}
Since this polynomial takes its values in $S$, it turns out that all
its coefficients belong to $S$. In particular, we have
\begin{eqnarray*}
\ad Y_i(X)=\partial_{r_i}\Phi(0)\in S\quad\mbox{and}\quad \big(\ad
Y_i\ad Y_j+\ad Y_j\ad
Y_i\big)X=2\partial_{r_i}\partial_{r_j}\Phi(0)\in S\,.
\end{eqnarray*}
The Jacobi identity can be read as
$\ad_U\ad_W-\ad_W\ad_U=\ad_{[U,W]}$, so that
\[
\big(\ad Y_i\ad Y_j+\ad Y_j\ad Y_i\big)X= 2\ad Y_i\ad Y_j
X+\ad[Y_j,Y_i]X.
\]
It follows that $\big(\ad Y_i\ad Y_j\big)X\in S$, and this proves
that $[\g',X]+[\g',[\g',X]]\subset S$. By induction, let us suppose
that
\[
\mathfrak{u}_{k-1}:=[\g',X]+[\g',[\g',X]]+\cdots
+[\underbrace{\g',[\g',\cdots,[\g'}_{\mbox{$(k-1)$
\small{times}}},X]\cdots]\subset S\,
\]
for some $k\geq 3$. In general we have
\begin{equation}\label{capodanno}
\partial_{r_{i_1}}\cdots\partial_{r_{i_k}}\Phi(0)=
\frac{1}{k!} \sum_{\sigma} \big(\ad Y_{j_{\sigma(1)}}\cdots \ad
Y_{j_{\sigma(k)}}\big)X\in S,
\end{equation}
where the sum runs on all permutations $\sigma$ of $k$ elements. By
the Jacobi identity $$\big(\ad Y_{j_{\sigma(1)}}\cdots\ad
Y_{j_{\sigma(k)}}\big)X-\big(\ad Y_{j_{\eta(1)}}\cdots\ad
Y_{j_{\eta(k)}}\big)X\in\mathfrak{u}_{k-1}$$ if
$\sigma\circ\eta^{-1}$ is a transposition. Then, by the inductive
assumption, we can iterate transpositions in $\big(\ad
Y_{j_{\sigma(1)}}\cdots\ad Y_{j_{\sigma(k)}}\big)X$ to write it as
$\big(\ad Y_{j_1}\cdots\ad Y_{j_k}\big)X+W_\sigma$ with $W_\sigma\in
S$. Then, from \eqref{capodanno} we get $\big(\ad Y_{j_1}\cdots\ad
Y_{j_k}\big)X\in S$. This shows that $\big(\ad Y_{j_1}\cdots\ad
Y_{j_k}\big)X\in S$, so that $\mathfrak{u}_k\subset S$, and this
proves the inclusion
$$ [\g',X]+
[\g',[\g',X]]+[\g',[\g',[\g',X]]]+\ldots\subset\span\{\Ad_{\exp
Y}X\mid Y\in\g'\}$$ Observing that the opposite inclusion trivially
holds, we are led to our claim.
\end{proof}

\subsection{Carnot groups}

A Carnot group $\G$ of step $s\geq 1$ is a connected, simply connected Lie
group whose Lie algebra $\g$ admits a step $s$ stratification: this
means that we can write
$$\g= V_1\oplus \cdots\oplus V_s $$
with $[V_j, V_1] = V_{j+1}$, $i\leq j\leq s$, $V_s\neq \{0\}$
and $V_{s+1}=\{0\}$. We keep the notation $n=\sum_i{\rm dim\,}V_i$
for the topological dimension of $\G$, and denote by
$$
Q:=\sum_{i=1}^s i\,{\rm dim\,}V_i
$$
the so-called \emph{homogeneous dimension} of $\G$. We denote by
$\delta_\lambda:\g\to\g$ the family of inhomogeneous dilations
defined by
$$
\delta_\lambda(\sum_{i=1}^s v_i):=\sum_{i=1}^s\lambda^i v_i
\qquad\lambda\geq 0
$$
where $X=\sum\limits_{i=1}^s v_i$ with $v_i\in V_i$, $1\leq i\leq
s$. The dilations $\delta_\lambda$ belong to $GL(\g)$ and are
uniquely determined by the homogeneity conditions
$$
\delta_\lambda X=\lambda^k X\quad\forall X\in V_k,\qquad 1\leq k\leq s.
$$

We denote by $m$ the dimension of $V_1$ and we fix an inner product
in $V_1$ and an orthonormal basis $X_1,\ldots,X_m$ of $V_1$. This
basis of $V_1$ induces the so-called Carnot-Caratheodory left
invariant distance $d$ in $\G$, defined as follows:
$$
d^2(x,y):=\inf\left\{\int_0^1\sum_{i=1}^m |a_i(t)|^2\,dt:\
\gamma(0)=x, \,\,\,\gamma(1)=y\right\},
$$
where the infimum is made among all Lipschitz curves
$\gamma:[0,1]\to \G$ such that $\gamma'(t)=\sum_1^m
a_i(t)(X_i)_{\gamma(t)}$ for a.e. $t\in [0,1]$ (the so-called
horizontal curves).

For Carnot groups, it is well known that the map $\exp: \g \to \G$
is a diffeomorphism, so any element $g\in \G$ can represented as
$\exp(X)$ for some unique $X\in\g$, and therefore uniquely written
in the form
\begin{equation}
\exp(\sum_{i=1}^s v_i),\qquad v_i\in V_i,\,\,1\leq i\leq s.
\end{equation}
This representation allows to define a family indexed by
$\lambda\geq 0$ of intrinsic dilations $\delta_\lambda:\G\to \G$, by
$$
\delta_\lambda\bigl(\exp(\sum_{i=1}^s v_i)\bigr):=
\exp\bigl(\sum_{i=1}^s \lambda^i v_i\bigr)\qquad (\text{i.e.
$\exp\circ\delta_\lambda=\delta_\lambda\circ\exp$.})
$$
We have kept the same notation $\delta_\lambda$ for both dilations
(in $\g$ and in $\G$) because no ambiguity will arise. Obviously,
$\delta_{\lambda}\circ\delta_{\eta}=\delta_{\lambda\eta}$, and the
Baker-Campbell-Hausdorff formula gives
$$
\delta_\lambda(xy)=\delta_\lambda (x)\delta_\lambda (y)
\qquad\forall x,\,y\in \G.
$$
Moreover, the Carnot-Caratheodory distance is well-behaved under these
dilations, namely
$$
d(\delta_\lambda x,\delta_\lambda y)=\lambda d(x,y)\qquad\forall x,\,y\in \G.
$$
Besides $\delta_\lambda\circ\exp=\exp\circ\delta_\lambda$, another
useful relation between dilations in $\G$ and dilations in $\g$ is
$\delta_\lambda X=(\delta_\lambda)_*X$, namely
\begin{equation}\label{valentino}
X(u\circ\delta_\lambda)(g)=(\delta_\lambda X)u(\delta_\lambda g)
\qquad\forall g\in\G,\,\,\lambda\geq 0.
\end{equation}
We have indeed
\begin{eqnarray*}
X(u\circ\delta_\lambda)(g)&=&\frac{d}{dt}
u\circ\delta_\lambda(g\exp(tX))\biggr\vert_{t=0}= \frac{d}{dt}
u(\delta_\lambda g\delta_\lambda\exp(tX))\biggr\vert_{t=0}\\&=&
\frac{d}{dt} u(\delta_\lambda g\exp(t\delta_\lambda
X))\biggr\vert_{t=0}=(\delta_\lambda X)u(\delta_\lambda g).
\end{eqnarray*}

\section{Measure-theoretic tools}

In this section we specify the notions of convergence used in this
paper (at the level of sets and of measures), and point out some
useful facts concerning Radon measures. The results quoted without
an explicit reference are all quite standard, and can be found for
instance in \cite{Ambook}, and those concerning Hausdorff measures
in metric spaces in \cite{Fed69} or \cite{AmbTil}.

\smallskip
{\bf Haar, Lebesgue and Hausdorff measures.} Carnot groups are
nilpotent and so unimodular, therefore the right and left Haar
measures coincide, up to constant multiples. We fix one of them
and denote it by $\vol_{\G}$.

We shall denote by ${\mathscr H}^k$ (resp. $\S^k$) the Hausdorff
(resp. spherical Hausdorff) $k$-dimensional measure; these measures
depend on the distance, and, unless otherwise stated, to build them
we will use the Carnot-Caratheodory distance in $\G$ and the
Euclidean distance in Euclidean spaces.

Using the left translation and scaling invariance of the
Carnot-Caratheodory distance one can easily check that the Haar
measures of $\G$ are a constant multiple of the spherical Hausdorff
measure $\S^Q$ and of ${\mathscr H}^Q$. In exponential coordinates,
all these measures are a constant multiple of the Lebesgue measure
$\L^n$ in $\R^n$, namely
$$
\vol_G\bigl(\{\exp(\sum_{i=1}^n x_iX_i):\ (x_1,\ldots,x_n)\in
A\}\bigr)=c \L^n(A)\qquad\text{for all Borel sets $A\subseteq\R^n$}
$$
for some constant $c$. Using this fact, one can easily prove that
\begin{equation}\label{scalQ}
\vol_\G(\delta_\lambda(A))=\lambda^Q\vol_\G(A)
\end{equation}
for all Borel sets $A\subseteq\G$.

The following implication will be useful: for $\mu$ nonnegative
Radon measure, $t>0$ and $B\subseteq\G$ Borel, we have
$$
\limsup_{r\downarrow 0}\frac{\mu(B_r(x))}{\omega_k r^k}\geq
t\,\,\,\forall x\in B\qquad \Longrightarrow\qquad \mu(B)\geq
t\S^k(B),
$$
where $\omega_k$ is the Lebesgue measure of the unit ball in $\R^k$
(it appears as a normalization constant in the definitions of
${\mathscr H}^k$ and $\S^k$, in order to ensure the identity
${\mathscr H}^k=\S^k=\L^k$ in $\R^k$). In particular we obtain that
\begin{equation}\label{cuzzola2}
\text{$\{x\in\G:\ \limsup_{r\downarrow
0}\frac{\mu(B_r(x))}{r^k}>0\}$ is $\sigma$-finite with respect to
$\S^k$.}
\end{equation}

\smallskip
{\bf Characteristic functions, convergence in measure.} For any set
$E$ we shall denote by $\uno_E$ the characteristic function of $E$
(1 on $E$, 0 on $\G\setminus E$); within the class of Borel sets of
$\G$, the convergence we consider is the so-called \emph{local
convergence in measure} (equivalent to the $L^1_{\rm loc}$
convergence of the characteristic functions), namely:
$$
E_h\to E\quad\Longleftrightarrow\quad \vol_\G\bigl(K\cap
[(E_h\setminus E)\cup (E\setminus E_h)]\bigr)=0\quad \text{for all
$K\subseteq \G$ compact.}
$$

\smallskip
{\bf Radon measures and their convergence.} The class ${\mathcal
M}(\G)$ of Radon measures in $\G$ coincides with the class of 0 order
distributions in $\G$, namely those distributions $T$ such that, for
any bounded open set $\Omega\subseteq \G$ there exists $C(\Omega)\in
[0,+\infty)$ satisfying
$$
|\langle T,g\rangle|\leq C(\Omega)\sup|g|\qquad\forall g\in
C^1_c(\Omega).
$$
These distributions can be uniquely extended to $C_c(\G)$, and their
action can be represented, thanks to Riesz theorem, through an
integral with respect to a $\sigma$-additive set function $\mu$
defined on bounded Borel sets. Thanks to this fact, the action of
these distributions can be extended even up to bounded Borel
functions with compact support. We will typically use both
viewpoints in this paper (for instance the first one plays a role in
the definition of distributional derivative, while the second one is
essential to obtain differentiation results). If $\mu$ is a
nonnegative Radon measure we shall denote
$$
{\rm supp\,}\mu:=\left\{x\in \G:\ \mu(B_r(x))>0\,\,\,\forall
r>0\right\}.
$$

The only convergence we use in ${\mathcal M}(\G)$ is the weak$^*$ one
induced by the duality with $C_c(\G)$, namely $\mu_h\to\mu$ if
$$
\lim_{h\to\infty}\int_\G u\,d\mu_h=\int_\G u\,d\mu\qquad\forall u\in
C_c(\G).
$$

\smallskip
{\bf Push-forward.} If $f:\G\to \G$ is a proper Borel map, then
$f^{-1}(B)$ is a bounded Borel set whenever $B$ is a bounded Borel
set. The push-forward measure $f_\sharp\mu$ is then defined by
$$
f_\sharp\mu(B):= \mu(f^{-1}(B)).
$$
In integral terms, this definition corresponds to
$$
\int_\G u\,d f_\sharp\mu:=\int_\G u\circ f\,d\mu
$$
whenever the integrals make sense (for instance $u$ Borel, bounded
and compactly supported).

\smallskip
{\bf Vector-valued Radon measures.} We will also consider
$\R^m$-valued Radon measures, representable as
$(\mu_1,\ldots,\mu_m)$ with $\mu_i\in{\mathcal M}(\G)$. The
\emph{total variation} of $|\mu|$ of an $\R^m$-valued measure $\mu$
is the smallest nonnegative measure $\nu$ defined on Borel sets of
$\G$ such that $\nu(B)\geq |\mu(B)|$ for all bounded Borel set $B$;
it can be explicitly defined by
$$
|\mu|(B):=\sup\left\{\sum_{i=1}^\infty |\mu(B_i)|:\ \text{$(B_i)$
Borel partition of $B$, $B_i$ bounded}\right\}.
$$

Push forward and convergence in ${\mathcal M}^m(\G)$ can be defined
componentwise. Useful relations between convergence and total
variation are:
\begin{equation}\label{savallo1}
\liminf_{n\to\infty}|\mu_n|(A)\geq|\mu|(A)\qquad\text{for all
$A\subseteq \G$ open,}
\end{equation}
\begin{equation}\label{savallo2}
\sup_{n\to\infty}|\mu_n|(K)<+\infty\qquad\text{for all $K\subseteq
\G$ compact,}
\end{equation}
whenever $\mu_n\to\mu$ in ${\mathcal M}^m(\G)$.

\smallskip
{\bf Asymptotically doubling measures.} A nonnegative Radon measure
$\mu$ in $\G$ is said to be \emph{asymptotically} doubling if
$$
\limsup_{r\downarrow 0}\frac{\mu(B_{2r}(x))}{\mu(B_r(x))}<+\infty
\qquad\text{for $\mu$-a.e. $x\in \G$.}
$$
For asymptotically doubling measures all the standard results of
Lebesgue differentiation theory hold: for instance, for any Borel
set $A$, $\mu$-a.e. point $x\in A$ is a \emph{density point} of $A$,
namely
$$
\lim_{r\downarrow 0}\frac{\mu(A\cap B_r(x))}{\mu(B_r(x))}=1.
$$
The same result holds for \emph{any} set $A$, provided we replace
$\mu$ by the outer measure $\mu^*$, defined for any $A\subseteq\G$
by
$$
\mu^*(A):=\inf\left\{\mu(B)\,:\,\text{$B$ Borel, $B\supseteq
A$}\right\}.
$$
It follows directly from the definition that $\mu^*$ is subadditive.
Moreover, let $(B_n)$ be a minimizing sequence and let $B$ the
intersection of all sets $B_n$: then $B$ is a Borel set, $B\supseteq
A$ and $\mu^*(A)=\mu(B)$. Furthermore, for all Borel sets $C$ we
have $\mu^*(A\cap C)=\mu(B\cap C)$ (if not, adding the strict
inequality $\mu^*(A\cap C)<\mu(B\cap C)$ to $\mu^*(A\setminus C)\leq
\mu(B\setminus C)$ would give a contradiction). Choosing $C=B_r(x)$, with
$x$ density point of $B$, we obtain
$$
\lim_{r\downarrow 0}\frac{\mu^*(A\cap B_r(x))}{\mu(B_r(x))}=
\lim_{r\downarrow 0}\frac{\mu(B\cap B_r(x))}{\mu(B_r(x))}=1.
$$
This proves that the set of points of $A$ that are not density
points is contained in a $\mu$--negligible Borel set. We will also
be using in the proof of Theorem~\ref{tpreiss} the fact that $\mu^*$
is \emph{countably} subadditive, namely
$\mu^*(A)\leq\sum_i\mu^*(A_i)$ for all sequences $(A_i)$ with
$A\subseteq\cup_iA_i$.

We recall the following result, proved in Theorem~2.8.17 of
\cite{Fed69}:

\begin{theorem}[Differentiation]\label{tfederer}
Assume that $\mu$ is asymptotically doubling and $\nu\in {\mathcal
M}(\G)$ is absolutely continuous with respect to $\mu$. Then the
limit
$$
f(x):=\lim_{r\downarrow 0}\frac{\nu(B_r(x))}{\mu(B_r(x))}
$$
exists and is finite for $\mu$-a.e. $x\in {\rm supp\,}\mu$.\\ In
addition, $f\in L^1_{\rm loc}(\mu)$ and $\nu=f\mu$, i.e.
$\nu(B)=\int_B f\,d\mu$ for all bounded Borel sets $B\subseteq\G$.
\end{theorem}

The proof given in \cite{Fed69} covers much more general situations;
the reader already acquainted with the theory of differentiation
with respect to doubling measures can easily realize that the
results extend to asymptotically doubling ones by consider the
\emph{localized} (in $\G\times (0,+\infty)$) maximal operators:
$$
M_{B,r}\nu(x):=\sup_{s\in (0,r)}\frac{\nu(B_s(x))}{\mu(B_s(x))},
\quad\qquad x\in B,
$$
where $\nu$ is any nonnegative Radon measure in $\G$. Thanks to the
asymptotic doubling property, one can find a family of Borel sets
$B_h\subseteq {\rm supp\,}\mu$ whose union covers $\G$, constants
$C_h\geq 1$ and radii $r_h>0$ such that $\mu(B_{3r}(x))\leq
C_h\mu(B_r(x))$ for $x\in B_h$ and $r\in (0,r_h)$. For the operators
$M_{B_h,r_h}$, the uniform doubling property on $B_h$ and a covering
lemma yield the weak $L^1$ estimate
$\mu(E\cap\{M_{B_h,r_h}\nu>t\})\leq t^{-1} C_h\nu(E)$ (for $E\subseteq B_h$
Borel, $t>0$). This leads to the differentiation result on all
$B_h$, and then $\mu$-a.e. on $\G$.

\section{Sets of locally finite perimeter}

In this section we recall a few useful facts about sets of finite
perimeter, considering also sets whose derivative along
non-horizontal directions is a measure.

\begin{definition}[Regular and invariant directions]
Let $f\in L^1_{\rm loc}(\G)$. \\We shall denote by ${\rm Reg}(f)$ the
vector subspace of $\g$ made by vectors $X$ such that $X f$ is
representable by a Radon measure. \\We shall denote by ${\rm
Inv}(f)$ the subspace of ${\rm Reg}(f)$ corresponding to the vector
fields $X$ such that $X f=0$, and by ${\rm Inv}_0(f)$ the subset
made by homogeneous directions, i.e.
$$
{\rm Inv}_0(f):={\rm Inv}(f)\cap\bigcup_{i=1}^s V_i.
$$
\end{definition}

Notice that, according to \eqref{righto},
$$f\circ
R_{\exp(tX)}=f\qquad\forall t\in\R,\,\, X\in {\rm Inv}(f).
$$

We will mostly consider regular and invariant directions of
characteristic functions, therefore we set
$$
{\rm Reg}(E):={\rm Reg}(\uno_E),\qquad {\rm Inv}(E):= {\rm
Inv}(\uno_E),\qquad {\rm Inv}_0(E):={\rm Inv}_0(\uno_E).
$$

We can now naturally define halfspaces by requiring invariance along
a codimension 1 space of directions, and monotonicity along the
remaining direction; if this direction is horizontal, we call these
sets \emph{vertical halfspaces}.

\begin{definition}[Vertical halfspaces]
We say that a Borel set $H\subseteq \G$ is a \emph{vertical
halfspace} if ${\rm Inv}_0(H)\supseteq\cup_2^s V_i$,
$V_1\cap {\rm Inv}_0(H)$ is a codimension one subspace
of $V_1$ and $X\uno_H\geq 0$ for some $X\in V_1$, with
$X\uno_H\neq 0$.
\end{definition}

Since
\begin{equation}\label{boissard}
{\rm span}\bigl({\rm Inv}_0(H)\bigr)=
\bigoplus_{i=1}^s V_i\cap {\rm Inv}_0(H),
\end{equation}
we can equivalently say that $H$ is an halfspace if
${\rm span}({\rm Inv}_0(H))$ is a codimension $1$ subspace
of $\g$, $V_1\cap {\rm span}({\rm Inv}_0(H))$ is a codimension
1 subspace of $V_1$ and $X\uno_H\geq 0$ for some $X\in V_1$:
indeed, \eqref{boissard} forces, whenever the
codimension is 1, all subspaces $V_i\cap {\rm Inv}_0(H)$ to coincide
with $V_i$, with just one exception.

Let us recall that $m$ denotes the dimension of $V_1$, and that
$X_1,\ldots,X_m$ is a given orthonormal basis of $V_1$. With this
notation, vertical halfspaces can be characterized as follows:

\begin{proposition}[Characterization of vertical
halfspaces]\label{pmaier} $H\subseteq\G$ is a vertical halfspace if
and only if there exist $c\in\R$ and a unit vector $\nu\in {\bf
S}^{m-1}$ such that $H=H_{c,\nu}$, where
\begin{equation}\label{defhcnu}
H_{c,\nu}:=\exp\bigl(\{\sum_{i=1}^m a_i X_i+\sum_{i=2}^s v_i:\
v_i\in V_i,\,\,a\in\R^m,\,\,\sum_{i=1}^m a_i\nu_i\leq c\}\bigr).
\end{equation}
\end{proposition}
\begin{proof}
Let us denote by $\nu\in {\bf S}^{m-1}$ the unique vector such that
the vector $Y=\sum_i\nu_iX_i$ is orthogonal to all invariant
directions in $V_1$.
Let us work in exponential coordinates, with the function
$$
(x_1,\ldots,x_n)\mapsto \exp(\sum_{i=1}^n x_iv_i),
$$
and let $\tilde{H}\subset\R^n$ be the set $H$ in these coordinates.
Here $(v_1,\ldots,v_n)$ is a basis of $\g$ compatible with the
stratification: this means that, if $m_i$ are the dimensions
of $V_i$, with $1\leq i\leq s$, $l_0=0$ and $l_i=\sum_1^i m_j$, then $v_{l_{i-1}+1},\ldots,v_{l_i}$
is a basis of $V_i$. By the Baker-Campbell-Hausdorff formula, in these coordinates
the vector fields $v_i$ correspond to $\partial_{x_i}$ for $l_{s-1}+1\leq i\leq l_s=n$,
and Theorem~\ref{tdp} gives that $\uno_{\tilde{H}}$ does not depend on
$x_{l_{s-1}+1},\ldots,x_n$. For $l_{s-2}+1\leq i\leq l_{s-1}$ the
vector fields $v_i-\partial_{x_i}$, still in these coordinates, are given
by the sum of polynomials multiplied by $\partial_{x_j}$, with $l_{s-1}+1\leq j\leq l_s$.
As a consequence $\partial_{x_i}\uno_{\tilde{H}}=0$ and we can apply Theorem~\ref{tdp}
again to obtain that $\uno_{\tilde{H}}$ does not depend on $x_{l_{s-2}+1},\ldots,x_{l_{s-1}}$
either. Continuing in this way we obtain that $\uno_{\tilde{H}}$
depends on $(x_1,\ldots,x_{m_1})$
only. Furthermore,
$\sum_i\xi_i\partial_{x_i}\uno_{\tilde{H}}$ is equal to $0$ if $\xi\perp\nu$,
and it is nonnegative if $\xi=\nu$. Then, a classical Euclidean argument
(it appears in De Giorgi's rectifiability
proof \cite{dg2}, see also the proof of this result in Theorem~3.59
of \cite{Ambook}) shows that $\uno_{\tilde{H}}$ depends on $\sum_1^m \nu_ix_i$
only, and it is a monotone function of this quantity. This immediately
gives \eqref{defhcnu}.
\end{proof}

\begin{remark} \label{rcambiaghi}{\rm
An analogous computation in exponential coordinates shows that ${\rm
Inv}(f)=\g$ if and only if $f$ is equivalent to a constant.}
\end{remark}

In the next proposition we point out useful stability properties of
${\rm Reg}(f)$ and ${\rm Inv}(f)$.

\begin{proposition} \label{pollo1} Let $f\in L^1_{\rm loc}(\G)$. Then
${\rm Reg}(f)$, ${\rm Inv}(f)$, ${\rm Inv}_0(f)$ are invariant under
left translations, and ${\rm Inv}_0(f)$ is invariant under intrinsic
dilations. Moreover:
\begin{itemize}
\item[(i)] ${\rm Inv}(f)$ is a Lie subalgebra of $\g$ and $[{\rm
Inv}_0(f),{\rm Inv}_0(f)]\subset{\rm Inv}_0(f)$;
\item[(ii)] If $X\in{\rm Inv}(f)$ and $k=\exp(X)$, then
$\Ad_k$ maps ${\rm Reg}(f)$ into ${\rm Reg}(f)$ and ${\rm Inv}(f)$
into ${\rm Inv}(f)$. More precisely
\begin{equation}\label{transfer}
\Ad_k(Y)f=(R_{k^{-1}})_\sharp Yf\qquad\forall Y\in {\rm Reg}(f).
\end{equation}
\end{itemize}
\end{proposition}
\begin{proof} The proof of the invariance is simple, so we omit it.

(i) We simply notice that for all $X,\,Y\in{\rm Inv}(f)$ we have
$$
\int_\G f [X,Y]g\,d\vol_\G= -\langle Xf,Yg\rangle+\langle
Yf,Xg\rangle=0 \qquad\forall g\in C^\infty_c(\G).
$$
The second stated property follows by the fact that
$[V_i,V_j]\subset V_{i+j}$.

(ii) Let $Y\in{\rm Reg}(f)$ and $Z=\Ad_k(Y)$. For $g\in
C^\infty_c(\G)$ and $k\in \G$ we have (taking into account the left
invariance of $Y$)
$$
Z g(x)=Y(g\circ C_k)(C_k^{-1}(x))=Y (g\circ R_{k^{-1}})(L_k\circ
C_{k^{-1}}(x))= Y (g\circ R_{k^{-1}})(R_k(x)).
$$
Therefore $(Zg)\circ R_{k^{-1}}= Y(g\circ R_{k^{-1}})$ and a change
of variables gives
$$
\int_\G f Zg\,d\vol_\G=\int_\G f\circ R_{k^{-1}} Y(g\circ
R_{k^{-1}})\,d\vol_\G.
$$
Now, if $k=\exp(X)$ with $X\in {\rm Inv}(f)$, we have $f\circ
R_{k^{-1}}=f$, and this gives \eqref{transfer}.
\end{proof}

\begin{remark}\label{rsiena}
{\rm Let $X\in {\rm Reg}(f)$ and assume that $Xf\geq 0$; then, combining
\eqref{eqknapp} with \eqref{transfer}, we obtain
$$
Xf+\sum_{i=1}^{s-1} \frac{t^i}{i!}\ad_Y^i(X)f\geq 0\qquad\forall t\in\R,\,\,
\forall Y\in{\rm Inv}(f).
$$
Since $t$ can be chosen arbitrarily large, this implies that
$$
\ad_Y^{s-1}(X)f\geq 0\qquad\forall Y\in{\rm Inv}(f).
$$
In particular, if $s$ is even, by applying the same inequality with $-Y$
in place of $Y$ we get
\begin{equation}\label{eqsiena}
\ad_Y^{s-1}(X)\in{\rm Inv}(f).
\end{equation}
}\end{remark}

\begin{definition}[Sets of locally finite perimeter]
The main object of investigation of this paper is the class of sets
of \emph{locally finite perimeter}, i.e. those Borel sets $E$ such that
$X\uno_E$ is a Radon measure for any $X\in V_1$.
\end{definition}

Still using the orthonormal basis of $V_1$, for $f\in L^1_{\rm
loc}(\G)$ with $X_if\in{\mathcal M}(\G)$ we can define the
$\R^m$-valued Radon measure
\begin{equation}\label{dchiE}
Df:=(X_1f,\ldots,X_mf).
\end{equation}
Two very basic properties that will play a role in the sequel are:
\begin{equation}\label{dnullo}
Df =0\qquad\Longrightarrow\qquad \text{$f$ is (equivalent to) a
constant}
\end{equation}
\begin{equation}\label{pcompa}
\text{$\sup_n\int_{\Omega}|f_n|\,d\vol_\G+|Df_n|(\Omega)<+\infty$\,\,
$\forall\Omega\Subset\G$}\quad\Longrightarrow\quad\text{$(f_n)$
relatively compact in $L^1_{\rm loc}$.}
\end{equation}
The proof of the first one can be obtained combining
Proposition~\ref{pollo1} (that gives that ${\rm Inv}(f)=\g$ with
Remark~\ref{rcambiaghi}). The second one has been proved in
\cite{garnieu}.

\begin{definition}[De Giorgi's reduced boundary]
Let $E\subseteq \G$ be a set of locally finite perimeter. We denote
by ${\mathscr F}E$ the set of points $x\in {\rm supp\,}|D\uno_E|$
where:
\begin{itemize}
\item[(i)] the limit $\nu_E(x)=(\nu_{E,1}(x),\ldots,\nu_{E,m}(x)):=\displaystyle{
\lim\limits_{r\downarrow 0}\frac{D\uno_E(B_r(x))}{|D\uno_E|(B_r(x))}}$ exists;
\item[(ii)] $|\nu_E(x)|=1$.
\end{itemize}
\end{definition}

The following result has been obtained in \cite{Amb01}.

\begin{theorem}\label{tambrosio} Let $E\subseteq \G$ be a set of locally finite
perimeter. Then $|D\uno_E|$ is asymptotically doubling, and more
precisely the following property holds: for $|D\uno_E|$-a.e.
$x\in\G$ there exists $\bar r(x)>0$ satisfying
\begin{equation}\label{cuzzola}
l_{\G} r^{Q-1}\leq |D\uno_E|(B_r(x))\leq L_{\G}r^{Q-1} \qquad\forall
r\in (0,\bar r(x)),
\end{equation}
with $l_{\G}$ and $L_{\G}$ depending on $\G$ only. As a consequence
$|D\uno_E|$ is concentrated on ${\mathscr F}E$, i.e.,
$|D\uno_E|(\G\setminus{\mathscr F}E)=0$.
\end{theorem}

Actually the result in \cite{Amb01} is valid in all Ahlfors
$Q$-regular metric spaces for which a Poincar\'e inequality holds
(in this context, obviously including all Lie groups, still the
measure $|D\uno_E|$ makes sense, see \cite{Mir}); \eqref{cuzzola}
also implies that the measure $|D\uno_E|$ can also be bounded from
above and below by the spherical Hausdorff measure ${\mathscr
S}^{Q-1}$, namely
\begin{equation}\label{cuzzola1}
\frac{l_{\G}}{\omega_{Q-1}}\S^{Q-1}(A\cap{\mathscr F}E)\leq
|D\uno_E|(A)\leq \frac{L_{\G}}{\omega_{Q-1}}\S^{Q-1}(A\cap{\mathscr
F}E)
\end{equation}
for all Borel sets $A\subseteq\G$ (since ${\mathscr H}^k\leq\S^k\leq
2^k{\mathscr H}^k$, similar inequalities hold with ${\mathscr
H}^{Q-1}$). In general doubling metric spaces, where no natural
dimension $Q$ exists, the asymptotic doubling property of
$|D\uno_E|$ and a suitable representation of it in terms of
Hausdorff measures have been obtained in \cite{Amb02}.

\section{Iterated tangents are halfspaces}

In this section we show that if we iterate sufficiently many times
the tangent operator we do get a vertical halfspace. Let us begin
with a precise definition of tangent set.

\begin{definition}[Tangent set]
Let $E\subseteq\G$ be a set of locally finite perimeter and $x\in
{\mathscr F}E$. We denote by ${\rm Tan}(E,x)$ all limit points, in
the topology of local convergence in measure, of the translated and
rescaled family of sets $\{\delta_{1/r}(x^{-1}E)\}_{r>0}$ as
$r\downarrow
0$.\\
If $F\in {\rm Tan}(E,x)$ we say that $F$ is tangent to $E$ at $x$.
We also set
$$
{\rm Tan}(E):=\bigcup_{x\in{\mathscr F}E}{\rm Tan}(E,x)
$$
\end{definition}

It is also useful to consider \emph{iterated} tangents; to this aim,
still for $x\in{\mathscr F}E$, we define ${\rm Tan}^1(E,x):={\rm
Tan}(E,x)$ and
$$
{\rm Tan}^{k+1}(E,x):=\bigcup\left\{{\rm Tan}(F):\ F\in {\rm Tan}^k(E,x)\right\}.
$$

The result we shall prove in this section is an intermediate step
towards Theorem~\ref{main}:

\begin{theorem}\label{main1}
Let $E\subseteq \G$ be a set with locally finite perimeter. Then, for
$|D\uno_E|$-a.e. $x\in \G$ we have (with the notation
\eqref{defhcnu})
$$
H_{0,\nu_E(x)}\in {\rm Tan}^k(E,x) \qquad\text{with}\qquad
k:=1+2(n-m).
$$
\end{theorem}

Notice that, by Theorem~\ref{tambrosio}, we need only to consider
points $x\in {\mathscr F}E$. Our starting point is the following
proposition, obtained in \cite{fsc}, showing that the tangent set at
points in the reduced boundary is always invariant along codimension
1 subspace of $V_1$, and monotone along the remaining horizontal
direction.

\begin{proposition}\label{fscc}
Let $E\subseteq \G$ be a set of locally finite perimeter. Then, for
all $\bar x\in {\mathscr F}E$ the following properties hold:
\begin{itemize}
\item[(i)] $0<\liminf\limits_{r\downarrow 0}|D\uno_E|(B_r(\bar x))/r^{Q-1}\leq
\limsup\limits_{r\downarrow 0}|D\uno_E|(B_r(\bar x))/r^{Q-1}<+\infty$;
\item[(ii)] ${\rm Tan}(E,\bar x)\neq\emptyset$ and, for all $F\in {\rm Tan}(E,\bar x)$,
we have that $e\in{\rm supp\,}|D\uno_F|$ and
$$
\nu_F=\nu_E(\bar x)\qquad\text{$|D\uno_F|$-a.e. in $\G$.}
$$
In particular $V_1\cap {\rm Inv}_0(F)$ coincides with the
codimension 1 subspace of $V_1$
$$
\left\{\sum_{i=1}^m a_iX_i:\
\sum_{i=1}^m a_i\nu_{E,i}(\bar x)=0\right\}
$$
and, setting, $X_x:=\sum_{i=1}^m \nu_{E,i}(\bar x)(X_i)_x\in\g$,
$X\uno_F$ is a nonnegative Radon measure.
\end{itemize}
\end{proposition}

In groups of step 2, in \cite{fsc} it is proved that constancy
of $\nu_E$ characterizes vertical subspaces. We provide here
a different proof of this fact, based on the properties of the
adjoint operator, and in particular on Remark~\ref{rsiena}.

\begin{proposition}\label{franchi}
Let $E\subset\G$ be a set with locally finite perimeter, and
assume that $\nu_E$ is (equivalent to) a constant. Then, if $\G$
is a step 2 group, $E$ is a vertical halfspace.
\end{proposition}
\begin{proof}
Let us denote by $\xi$ the constant value of $\nu_E$, and set
$X:=\sum_i\xi_iX_i$. Then $X\uno_E\geq 0$ and ${\rm Inv}(E)$
contains all vectors $Y=\sum_i\eta_iX_i$ with $\eta\in\R^m$
perpendicular to $\xi$. From \eqref{eqsiena} we get $[Y,X]\uno_E=0$
for any $Y\in {\rm Inv}(E)\cap V_1$, and since these commutators,
together with the commutators $\{[Y_1,Y_2]:\ Y_i\in{\rm Inv}(E)\cap
V_1\}$, span the whole of $V_2$, the proof is achieved.
\end{proof}

\begin{remark}\label{rloca}{\rm
The following simple example, that we learned from F. Serra Cassano,
shows that the sign condition is essential for the validity of the
classification result, even in the first Heisenberg group $\H^1$.
Choosing exponential coordinates $(x,y,t)$, and the vector fields
$X_1:=\partial_x+2y\partial_t$ and $X_2:=\partial_y-2x\partial t$,
the function
$$
f(x,y,t):=g(t+2xy)
$$
(with $g$ smooth) satisfies $X_1f=4yg'(t+2xy)$ and $X_2f=0$. Therefore the
sets $E_t:=\{f<t\}$ are $X_2$-invariant and are not halfspaces. The same example
can be used to show that there is no local version of Proposition~\ref{franchi},
because the sets $E_t$ locally may satisfy $X_1\uno_{E_t}\geq 0$ or
$X_1\uno_{E_t}\leq 0$ (depending on the sign of $g'$ and $y$), but
are not locally halfspaces.

The non-locality appears also in our argument: indeed, the proof of \eqref{eqsiena} depends
on the sign condition of $\Ad_{\exp(tX_2)}(X_1)\uno_E$ with $t$ arbitrarily large, and this
is the right translate, by $\exp(tX_2)$, of $X_1\uno_E$. The proof given in \cite{fsc}
depends, instead, on the possibility of joining two different points in $\H^1$ by following
integral lines of $X_2$ in both directions, and integral lines of $X_1$ in just
one direction: an inspection of the proof reveals that these paths can not be confined
in a bounded region, even if the initial and final point are confined within a
small region. In this sense, Proposition~\ref{franchi} could be considered as
a kind of Liouville
theorem.
}\end{remark}

Let $f\in L^1_{\rm loc}(\G)$ and $X\in\g$; then, for all $r>0$ we
have the identity
\begin{equation}\label{inhomo}
\delta_{1/r}X (f\circ \delta_r)= r^{-Q} (\delta_{1/r})_\sharp(Xf)
\end{equation}
in the sense of distributions. Indeed, writing in brief
$X_r:=\delta_{1/r}X$, if $g\in C^\infty_c(\G)$, from
\eqref{valentino} we get $X_r(g\circ\delta_r)= (Xg)\circ\delta_r$;
as a consequence \eqref{scalQ} gives
\begin{eqnarray}\label{scalo}
\langle X_r (f\circ\delta_r), g\rangle&=& -\int_\G (f\circ\delta_r)
X_r g\,d\vol_\G=-r^{-Q}\int_\G f (X_r g)\circ \delta_{1/r}\,d\vol_\G
\\&=&-r^{-Q}\int_\G f X (g\circ\delta_{1/r})\,d\vol_Q=\langle r^{-Q}
(\delta_{1/r})_\sharp(Xf),g\rangle.\nonumber
\end{eqnarray}

The first crucial lemma shows that if $X\in {\rm Reg}(E)$ belongs to
$\oplus_2^s V_i$, then the tangents to $E$ at $|D\uno_E|$-a.e. $x$
are invariant under $Y$, where $Y$ is the ``higher degree
part'' of $X$ induced by the stratification of $\g$. The underlying
reason for this fact is that the intrinsic dilations behave quite
differently in the $X$ direction and in the horizontal direction.

\begin{lemma}\label{provafis}
Let $F$ be a set with locally finite perimeter, $X\in {\rm Reg}(F)$,
$\mu=X\uno_F$ and assume that $X=\sum_{i=2}^l v_i$ with $v_i\in V_i$
and $l\leq s$. Then, for $|D\uno_F|$-a.e. $x$, $v_l\in{\rm
Inv}_0(L)$ for all $L\in{\rm Tan}(F,x)$.
\end{lemma}
\begin{proof}
>From \eqref{cuzzola2} we know that the set $N$ of points $x$ such
that $\limsup_{r\downarrow 0} r^{2-Q}|\mu|(B_r(x))$ is positive is
$\sigma$-finite with respect to $\S^{Q-2}$, and therefore
$\S^{Q-1}$-negligible and $|D\uno_F|$-negligible (recall
\eqref{cuzzola1}). We will prove that the statement holds at any
$x\in ({\mathscr F}F)\setminus N$ and we shall assume, up to a left
translation, that $x=e$. Given any $g\in C^1_c(\G)$, let $R$ be such
that ${\rm supp}(g)\subsetneq B_R(e)$; \eqref{inhomo} with
$f=\uno_F$ gives
$$
\int_\G \uno_{\delta_{1/r}F} X_rg\,d\vol_\G=r^{l-Q}\int_\G
g\circ\delta_{1/r}\,d\mu
$$
with $X_r:=r^l\delta_{1/r}X$, so that $X_r\to v_l$ as $r\downarrow 0$.
Now, notice that $l\geq 2$, and that
the right hand side can be bounded with
$$
\sup|g| r^{l-Q}|\mu|(B_{Rr}(e))= O(r^{l-Q})o(r^{Q-2})=o(1).
$$
So, passing to the limit as $r\downarrow 0$ along a suitable
sequence, we obtain that $v_l\uno_L=0$ for all $L\in{\rm Tan}(F,e)$.
\end{proof}

The invariance of ${\rm Inv}_0$ under left translations and scaling
shows that ${\rm Inv}_0(F)$ contains ${\rm Inv}_0(E)$ for all $F\in
{\rm Tan}(E)$. Let us define codimension of ${\rm Inv}_0(E)$ in $\g$
as the codimension of its linear span;
we know that this codimension is at least 1 (because the
codimension within $V_1$ is 1) for \emph{all} tangent sets, and it is
equal to 1 precisely for vertical halfspaces, thanks to
Proposition~\ref{pmaier}.

The second crucial lemma shows that, when the codimension of ${\rm
Inv}_0(E)$ in $\g$ is at least 2, a double tangent \emph{strictly}
increases, at $|D\uno_E|$-a.e. point, the set ${\rm Inv}_0(E)$. The
strategy is to find first a tangent set $F$ with ${\rm
Reg}(F)\supsetneq {\rm span}({\rm Inv}_0(E))$ (this is based on the
geometric Proposition~\ref{pcrucial} and Proposition~\ref{pollo1})
and then on the application of the previous lemma, which turns a
regular direction of $F$ into an invariant homogeneous direction of
a tangent to $F$.

\begin{lemma}[Improvement of ${\rm Inv}_0(E)$]\label{limpro2}
Let $E\subseteq\G$ be a set of locally finite perimeter and assume
that
$$
{\rm dim}\bigl({\rm span}({\rm Inv}_0(E))\bigr)\leq n-2.
$$
Then, for all $\bar x\in {\mathscr F}E$, ${\rm
Inv}_0(L)\supsetneq {\rm Inv}_0(E)$ for some $L\in {\rm
Tan}^2(E,\bar x)$.
\end{lemma}
\begin{proof} {\bf (Step 1)} We show first the existence of $Z\in\g\setminus
[{\rm span}({\rm Inv}_0(E))+V_1]$ such that $Z\in {\rm Reg}(F)$ for
all $F\in {\rm Tan}(E,\bar x)$. To this aim, we apply
Proposition~\ref{pcrucial} with $\g':={\rm span}({\rm Inv}_0(E))$
(recall that, by Proposition~\ref{pollo1}(i), $\g'$ is a Lie
algebra) and $X:=\sum_1^m\nu_{E,i}(\bar x)X_i$ to obtain $Y\in\g'$
such that
$$
Z:=\Ad_{\exp(Y)}(X) \notin {\rm span}({\rm Inv}_0(E))\oplus \{\R
X\}= {\rm span}({\rm Inv}_0(E))+V_1.$$ Then, since ${\rm Inv}_0(F)$
contains ${\rm Inv}_0(E)$ for all $F\in{\rm Tan}(F,\bar x)$, we have
that $Y\in{\rm Inv}(F)$, therefore Proposition~\ref{pollo1}(ii)
shows that $Z\in {\rm Reg}(F)$ for all $F\in {\rm Tan}(E,\bar x)$.

{\bf (Step 2)} Now, let $F\in {\rm Tan}(E,\bar x)$, $Z\notin {\rm
span}({\rm Inv}_0(E))+V_1$ given by the previous step, and set
$\mu=Z\uno_F$. Possibly removing from $Z$ its horizontal component
we can write $Z=v_{i_1}+\cdots+v_{i_l}$ with $i_j\geq 2$ and
$v_{i_j}\in V_{i_j}$. Then, $v_{i_k}\notin{\rm Inv}_0(E)$ for at
least one $k\in\{1,\ldots,l\}$, and let us choose the largest one
with this property. Then, setting $Z'=v_{i_1}+\cdots+v_{i_k}$, since
$v_{i_j}\in{\rm Inv}_0(E)\subseteq{\rm Inv}_0(F)$ for all $k<j\leq
l$, we still have $Z'\uno_F=\mu$. By Lemma~\ref{provafis} we can
find $L\in{\rm Tan}(F)$ with $v_{i_k}\uno_L=0$, i.e. $v_{i_k}\in
{\rm Inv}_0(L)$. Since $v_{i_k}\notin{\rm Inv}_0(E)$, we have proved
that ${\rm Inv}_0(L)$ strictly contains ${\rm Inv}_0(E)$.
\end{proof}

{\bf Proof of Theorem~\ref{main1}.} Recall that $m={\rm dim}(V_1)$.
Sets in ${\rm Tan}(E,\bar x)$ are invariant, thanks to
Proposition~\ref{fscc}, in at least $m-1$ directions. Let us define
$$
i_k:=\max\left\{{\rm dim}\bigl({\rm span}({\rm Inv}_0(F))\bigr):\
F\in {\rm Tan}^k(E,\bar x)\right\}.
$$
Then $i_1\geq m-1$ and we proved in Lemma~\ref{limpro2} that
$i_{k+2}>i_k$ as long as there exists $F\in {\rm Tan}^k(E,\bar x)$
with ${\rm dim}\bigl({\span}({\rm Inv}_0(F))\bigr)\leq n-2$.
By iterating $k$ times, with
$k\leq 2(n-m)$, the tangent operator we find $F\in {\rm Tan}^k(E,\bar x)$
with ${\rm dim}\bigl({\span}({\rm Inv}_0(F))\bigr)\geq n-1$.

We know from Proposition~\ref{fscc} that $e\in{\rm supp\,}|D\uno_F|$, that
the codimension of ${\rm Inv}_0(F)$ is exactly 1, and precisely that
$$
V_1\cap {\rm Inv}_0(F)=\left\{\sum_{i=1}^m a_iX_i:\
\sum_{i=1}^m a_i\nu_{E,i}(\bar x)=0\right\}
$$
and that $\sum_i\nu_{E,i}(\bar x)X_i\uno_F\geq 0$.
Therefore Proposition~\ref{pmaier} gives $F=H_{0,\nu_E(\bar x)}$.

\section{Iterated tangents are tangent}\label{sitera}

In this section we complete the proof of Theorem~\ref{main}. Taking
into account the statement of Theorem~\ref{main1}, we need only to
prove the following result.

\begin{theorem}\label{main2}
Let $E\subseteq \G$ be a set with locally finite perimeter. Then, for
$|D\uno_E|$-a.e. $x\in \G$ we have
$$
\bigcup_{k=2}^\infty {\rm Tan}^k(E,x)\subseteq{\rm Tan}(E,x).
$$
\end{theorem}

In turn, this result follows by an analogous one involving tangents
to measures, proved in \cite{Preiss} in the Euclidean case; we just
adapt the argument to Carnot groups and to vector-valued measures.
In the sequel we shall denote by $I_{x,r}(y):=\delta_{1/r}(x^{-1}y)$
the composition $\delta_{1/r}\circ L_{x^{-1}}$.

We say that a measure $\mu\in{\mathcal M}^m(\G)$ is
\emph{asymptotically $q$-regular} if
\begin{equation}\label{capodanno5}
0<\liminf_{r\downarrow 0}\frac{|\mu|(B_r(x))}{r^q} \leq
\limsup_{r\downarrow 0}\frac{|\mu|(B_r(x))}{r^q}<+\infty
\qquad\text{for $|\mu|$-a.e. $x\in\G$.}
\end{equation}
Notice that asymptotically $q$-regular measures are asymptotically
doubling, and that the perimeter measure $|D\uno_E|$ is
asymptotically $(Q-1)$-regular, thanks to Theorem~\ref{tambrosio}.

\begin{definition}[Tangents to a measure]
Let $\mu\in {\mathcal M}^m(\G)$ be asymptotically $q$-regular. We
shall denote by ${\rm Tan}(\mu,x)$ the family of all measures
$\nu\in {\mathcal M}^m(\G)$ that are weak$^*$ limit points as
$r\downarrow 0$ of the family of measures
$r^{-q}(I_{x,r})_\sharp\mu$.
\end{definition}

\begin{theorem}\label{tpreiss}
Let $\mu\in{\mathcal M}^m(\G)$ be asymptotically $q$-regular. Then,
for $|\mu|$-a.e. $x$, the following property holds:
$$
{\rm Tan}(\nu,y)\subseteq {\rm Tan}(\mu,x)\qquad \forall \nu\in{\rm
Tan}(\mu,x),\,\,y\in{\rm supp\,}|\nu|.
$$
\end{theorem}

The connection between Theorem~\ref{main2} and Theorem~\ref{tpreiss}
rests on the following observation:
\begin{equation}\label{equivalence}
L\in {\rm Tan}(F,x)\qquad\Longleftrightarrow\qquad D\uno_L\in{\rm
Tan}(D\uno_F,x)\setminus\{0\}
\end{equation}
for all $x\in {\mathscr F}F$. The implication $\Rightarrow$ in
\eqref{equivalence} is easy, because a simple scaling argument gives
\begin{equation}\label{equivalencebis}
L=\lim_{i\to\infty}\delta_{1/r_i}(x^{-1}F)\qquad\Longrightarrow\qquad
D\uno_L=\lim_{i\to\infty} r_i^{1-Q} (I_{x,r_i})_\sharp D\uno_F.
\end{equation}
Therefore $L\in {\rm Tan}(F,x)$ implies $D\uno_L\in{\rm
Tan}(D\uno_F,x)$; clearly $D\uno_L\neq 0$ because $x\in {\mathscr
F}F$.

Now we prove the harder implication $\Leftarrow$ in
\eqref{equivalence}: assume, up to a left translation, that $x=e$,
and that $D\uno_L\neq 0$ is the weak$^*$ limit of
$r_i^{1-Q}(I_{e,r_i})_\sharp D\uno_F$, with $r_i\downarrow 0$; now,
set $F_i:=\delta_{1/r_i}F$, so that $D\uno_{F_i}=r_i^{1-Q}
(I_{e,r_i})_\sharp D\uno_F$, and by the compactness properties of
sets of finite perimeter (see \eqref{pcompa}) assume with no loss of
generality that $F_i\to L'$ locally in measure, so that $L'\in{\rm
Tan}(F,e)$. Then $r_i^{1-Q}(I_{e,r_i})_\sharp D\uno_F= D\uno_{F_i}$
weakly$^*$ converge to $D\uno_{L'}$: indeed, the convergence in the
sense of distributions is obvious, and since the total variations
are locally uniformly bounded, we have weak$^*$ convergence as well.
It follows that $D\uno_L=D\uno_{L'}$. Since $\uno_L-\uno_{L'}$ has
zero horizontal distributional derivative, by \eqref{dnullo} it must
be (equivalent to) a constant; this can happen only when either
$L=L'$ or $L=\G\setminus L'$; but the second possibility is ruled
out because it would imply that $D\uno_L=-D\uno_{L'}$ and that
$D\uno_L=0$. This proves that $L=L'\in{\rm Tan}(F,e)$.

\smallskip
\begin{proof} {\bf (of Theorem~\ref{main2})} At any point $x\in{\mathscr
F}E$ where the property stated in Theorem~\ref{tpreiss} holds with
$\mu=D\uno_E$ we may consider any $F\in{\rm Tan}(E,x)$ and $L\in
{\rm Tan}(F,y)$ for some $y\in {\mathscr F}F$; then, by
\eqref{equivalence} we know that $D\uno_F\in {\rm Tan}(D\uno_E,x)$
and $D\uno_L\in {\rm Tan}(D\uno_F,y)\setminus\{0\}$; as a
consequence, Theorem~\ref{tpreiss} gives $D\uno_L\in{\rm
Tan}(D\uno_E,x)\setminus\{0\}$, hence \eqref{equivalence} again
gives that $L\in {\rm Tan}(E,x)$. This proves that ${\rm
Tan}^2(E,x)\subseteq {\rm Tan}(E,x)$, and therefore ${\rm
Tan}^3(E,x)\subseteq {\rm Tan}^2(E,x)$, and so on.
\end{proof}

The rest of this section is devoted to the proof of
Theorem~\ref{tpreiss}. We will follow with minor variants (because
we are dealing with vector-valued measures) the proof given in
Mattila's book \cite{mattila}. Before proceeding to the proof of
Theorem~\ref{tpreiss} we state a simple lemma.

\begin{lemma}\label{ldensity}
Assume that $A\subset\G$ and $a\in A$ is a density point for $A$
relative to $|\mu|^*$, i.e.
\begin{equation}
\lim_{r\downarrow 0}\dfrac{|\mu|^*(B_r(a)\cap A)}{|\mu|(B_r(a))}=1.
\end{equation}
If, for some $r_i\downarrow 0$ and $\lambda_i\geq 0$ the measures
$\lambda_i(I_{a,r_i})_\sharp\mu$ weakly$^*$ converge to $\nu$, then
$$
\lim_{i\to\infty}\frac{d(a\delta_{r_i}y,A)}{r_i}=0 \qquad\forall
y\in {\rm supp\,}|\nu|.
$$
\end{lemma}
\begin{proof}
Let $\tau:=d(y,e)$ and let us argue by contradiction. If the
statement were false, $\tau$ would be positive and there would exist
$\eps\in (0,\tau)$ such that $d(a\delta_{r_i}y,A)>\eps r_i$ for
infinitely many values of $i$. Possibly extracting a subsequence,
let us assume that this happens for all $i$: we know that
\begin{equation}\label{mattila1}
B_{\eps r_i}(a\delta_{r_i}y)\subseteq \G\setminus A
\end{equation}
and since $\eps<\tau$ we have
\begin{equation}\label{mattila2}
B_{\eps r_i}(a\delta_{r_i}y)\subseteq B_{\tau r_i}(a\delta_{r_i}y)
\subseteq B_{2\tau r_i}(a).
\end{equation}
Now use, in this order, the definition of density point,
\eqref{mattila1}, \eqref{mattila2} and \eqref{savallo1} to get
\begin{eqnarray*}
1&=&\lim_{i\to\infty}\dfrac{|\mu|^*(B_{2\tau r_i}(a))\cap
A)}{|\mu|(B_{2\tau r_i}(a))}
\leq\limsup_{i\to\infty}\dfrac{|\mu|(B_{2\tau r_i}(a)\setminus
B_{\eps r_i}(a\delta_{r_i}y))}{|\mu|(B_{2\tau r_i}(a))}\\
&=&\limsup_{i\to\infty}\dfrac{|\mu|(B_{2\tau r_i}(a)))-|\mu|
(B_{\eps r_i}(a\delta_{r_i}y))}{|\mu|(B_{2\tau r_i}(a))}
=1-\liminf_{i\to\infty}\dfrac{|\mu| (B_{\eps r_i}(a\delta_{r_i}y))}{|\mu|(B_{2\tau r_i}(a))}\\
&=&1-\liminf_{i\to\infty}\dfrac{(I_{a,r_i})_\sharp|\mu|
(B_\eps(y))}{(I_{a,r_i})_\sharp|\mu|(B_{2\tau}(e))}\leq
1-\dfrac{\liminf_{i\to\infty}|\lambda_i(I_{a,r_i})_\sharp\mu|
(B_\eps(y))}
{\limsup_{i\to\infty}|\lambda_i(I_{a,r_i})_\sharp\mu|(B_{2\tau}(e))}\\
&\leq&1-
\dfrac{|\nu|(B_\eps(y))}{\limsup_{i\to\infty}|\lambda_i(I_{a,r_i})_\sharp\mu|(B_{2\tau}(e))}.
\end{eqnarray*}
But, $|\nu|(B_\eps(y))>0$ because $y\in {\rm supp\,}|\nu|$,
and the $\limsup$ is finite by \eqref{savallo2}. This contradiction
concludes the proof of the lemma.
\end{proof}

{\bf Proof of Theorem~\ref{tpreiss}.} For $\nu,\,\nu'\in {\mathcal
M}^m(\G)$, define
$$
d_R(\nu,\nu'):=\sup\left\{\int_{\G}\phi\,d\nu-\int_{\G}\phi\,d\nu':\
\phi\in {\mathcal D}_R\right\},
$$
where
$${\mathcal D}_R:=\left\{\phi\in C_c(B_R(e)):\
\sup|\phi|\leq 1\quad\text{and}\quad |\phi(x)-\phi(y)|\leq d(x,y)
\,\,\forall x,\,y\in\G\right\}. $$ It is well known, and easy to
check, that $d_R$ induces the weak$^*$ convergence in all bounded
sets of ${\mathcal M}^m(B_R(e))$. We define a distance $\bar d$ in
${\mathcal M}^m(\G)$ by
$$
\bar d(\mu,\nu):=\sum_{R=1}^\infty
2^{-R}\min\bigl\{1,d_R(\mu,\nu)\bigr\}.
$$
Let $x$ be a point where the limsup in
\eqref{capodanno5} is finite; now we
check that, for all infinitesimal sequences $(r_i)\subset
(0,+\infty)$, we have
\begin{equation}\label{capodanno1}
\nu={\rm weak^*}-\lim_{i\to\infty}r_i^{-q}
(I_{x,r_i})_\sharp\mu\qquad\Longleftrightarrow
\qquad\lim_{i\to\infty}\bar d(\nu,r_i^{-q}(I_{x,r_i})_\sharp\mu)=0.
\end{equation}
The implication $\Rightarrow$ is obvious, because $\bar
d$-convergence is equivalent to $d_R$-convergence for all $R$, and
all weakly$^*$-convergent sequences are locally uniformly bounded
(see \eqref{savallo2}). The implication $\Leftarrow$ is analogous,
but it depends on our choice of $x$, which ensures the property
$$
\sup_{i\in\N}r_i^{-q}|(I_{x,r_i})_\sharp\mu|(B_R(e))= \sup_{i\in\N}
r_i^{-q}|\mu|(B_{Rr_i}(x))\leq R^q\limsup_{r\downarrow 0}\frac{|\mu|(B_r(x))}{r^q}
<+\infty.
$$
This property ensures that $r_i^{-q}(I_{x,r_i})_\sharp\mu$ is
bounded in all ${\mathcal M}^m(B_R(e))$ for all $R>0$, and enables
to pass from $d_R$-convergence to weak$^*$ convergence in all balls
$B_R(e)$.

Thanks to the equivalence stated in \eqref{capodanno1}, by a
diagonal argument it suffices to prove that, for $|\mu|$-a.e. $x$,
the following property holds: for all $\nu\in {\rm Tan}(\mu,x)$,
$y\in {\rm supp\,}|\nu|$ and $r>0$ we have
$r^{-q}(I_{y,r})_\sharp\nu\in {\rm Tan}(\mu,x)$. But since the
operation $\sigma\mapsto r^{-q}(I_{e,r})_\sharp\sigma$ is easily
seen to map ${\rm Tan}(\mu,x)$ into ${\rm Tan}(\mu,x)$, and
$I_{y,r}=I_{e,r}\circ I_{y,1}$, we need just to show that:

\noindent (*) for $|\mu|$-a.e. $x$ the following property holds: for
all $\nu\in {\rm Tan}(\mu,x)$ and all $y\in {\rm supp\,}|\nu|$, we
have $(I_{y,1})_\sharp\nu\in {\rm Tan}(\mu,x)$.

Heuristically, this property holds at ``Lebesgue'' points of the
multivalued map $x\mapsto {\rm Tan}(\mu,x)$, thanks to the identity
\begin{equation}\label{chainI}
I_{\delta_{1/r}(x^{-1}y),1}\circ I_{x,r}=I_{y,r}.
\end{equation}
Indeed, this identity implies that tangents to $\mu$ at $x$ on the
scale $r$ are close to tangents to $\mu$ at $y$ on the scale $r$
when $d(x,y)\ll r$.

Let us consider the set $R$ of points where the property (*) fails:
for all $x\in R$ there exist a measure $\nu\in {\rm Tan}(\mu,x)$ and
a point $y\in {\rm supp\,}|\nu|$ such that
$(I_{y,1})_\sharp\nu\notin {\rm Tan}(\mu,x)$. This implies, thanks
to the implication $\Leftarrow$ in \eqref{capodanno1}, the existence
of integers $z,\,k\geq 1$ such that the measure
$(I_{y,1})_\sharp\nu$ is $1/k$ far (relative to $\bar d$) from the
set $r^{-q}(I_{x,r})_\sharp\mu:\ r\in (0,1/z) \}$. Set
$$
A_{z,k}:=\left\{x\in\G:\ \exists \nu\in {\rm Tan}(\mu,x),\,\,
\exists y\in {\rm supp\,}|\nu|\,\,\text{such that}\right.$$
$$\left.\overline{d}((I_{y,1})_\sharp\nu,r^{-q}
(I_{x,r})_\sharp\mu)>1/k,\,\,\,\forall r\in (0,1/z)\right\}.$$ Since
$R$ is contained in the union of these sets, to conclude the proof
it suffices to show that $|\mu|^*(A_{z,k})=0$ for any $z,\,k\geq 1$.

Suppose by contradiction $|\mu|^*(A_{z,k})\neq0$ for some $z,\,k\geq
1$ and let us fix these two parameters; it is not difficult to check
that we can cover the space ${\mathcal M}^m(\G)$ with a family
$\{B_l\}$ of sets satisfying
\begin{equation}\label{capodanno4}
\bar d(\nu,\nu')<\frac{1}{2k}\qquad\forall \nu,\,\nu'\in B_l.
\end{equation}
Let us now consider the sets
$$
A_{z,k,l}:=\left\{x\in\G:\ \exists \nu\in {\rm Tan}(\mu,x),\,\,
\exists y\in {\rm supp\,}|\nu|\,\,\text{such
that}\,\,(I_{y,1})_\sharp\nu\in B_l,\right.$$
$$\left.
\overline{d}((I_{y,1})_\sharp\nu,r^{-q}(I_{x,r})_\sharp\mu)>1/k,\,\,\,
\forall r\in (0,1/z)\right\}.$$ Since $\cup_l A_{z,k,l}$ contains
$A_{z,k}$ and $|\mu|^*$ is countably subadditive, at least one of
these sets satisfies $|\mu|^*(A_{z,k,l})>0$. Let us fix $l$ with
this property, and let us denote $A_{z,k,l}$ by $A$.

Since $|\mu|^*(A)>0$ and $|\mu|$ is asymptotically doubling, we can
find $a\in A$ which is a density point of $A$ relative to $|\mu|^*$.
>From now on also the point $a$ will be fixed, and so an associated
measure $\nu_a\in {\rm Tan}(\mu,a)$, a point $y_a\in {\rm
supp\,}|\nu_a|$ satisfying $(I_{y_a,1})_\sharp\nu_a\in B_l$ and
\begin{equation}\label{eccoa}
\overline{d}((I_{y_a,1})_\sharp\nu_a,r^{-q}(I_{a,r})_\sharp\mu)>\frac{1}{k},\,\,\,
\forall r\in (0,1/m).
\end{equation}
We can also write $\nu_a=\lim_{i\to\infty}
r_i^{-q}(I_{a,r_i})_\sharp\mu$, for suitable $r_i\downarrow 0$, and
clearly \eqref{eccoa} implies that $y_a\neq e$.

Let us consider the points $a\cdot\delta_{r_i}y_a$ and their
distance from $A$ and take $a_i\in A$ such that $\dist
(a\delta_{r_i}y_a,a_i)\leq\dist (a\delta_{r_i}y_a,A)+r_i/i$.
Lemma~\ref{ldensity} yields that $\dist
(a\delta_{r_i}y_a,a_i)=o(r_i)$ as $i\to\infty$, and so
$\delta_{1/r_i}(a^{-1}a_i)\to y_a$. Now, \eqref{chainI} shows that
$I_{\delta_{1/r_i}(a^{-1}a_i), 1}\circ I_{a,r_i}=I_{a_i,r_i}$, so
that
\begin{eqnarray*}
\lim_{i\to\infty}r_i^{-q}(I_{a_i,r_i})_\sharp\mu&=&\lim_{i\to\infty}
r_i^{-q}(I_{\delta_{1/r_i}(a^{-1}a_i), 1})_\sharp(I_{a,r_i})_\sharp\mu\\
&=&\lim_{i\to\infty}(I_{\delta_{1/r_i}(a^{-1}a_i),
1})_\sharp\bigl(r_i^{-q}(I_{a,r_i})_\sharp\mu\bigr) =
(I_{y_a,1})_\sharp\nu_a.\end{eqnarray*} So, we can fix $i$
sufficiently large such that $r_i<1/z$ and
\begin{equation}\label{II}\overline{d}(r_i^{-q}(I_{a_i,r_i})_\sharp\mu,
(I_{y_a,1})_\sharp\nu_a)<\dfrac{1}{2k}.\end{equation} Since $a_i\in
A= A_{z,k,l}$, we can find a measure $\nu'\in {\rm Tan}(\mu,a_i)$
and a point $y'\in {\rm supp\,}|\nu'|$ with
$(I_{y',1})_\sharp\nu'\in B_l$ such that
$$
\dfrac{1}{k}<\overline{d}(r_i^{-q}(I_{a_i,r_i})_\sharp\mu,(I_{y',1})_\sharp\nu').$$
By applying the triangle inequality we obtain
$$
\frac{1}{k}<\overline{d}(r_i^{-q}(I_{a_i,r_i})_\sharp\mu,
(I_{y_a,1})_\sharp\nu_a)
+\overline{d}((I_{y_a,1})_\sharp\nu_a,(I_{y',1})_\sharp\nu')<
\frac{1}{2k}+\frac{1}{2k},
$$
where we used \eqref{II} and our choice \eqref{capodanno4} of $B_l$.
The contradiction ends the proof of the theorem.

\section{The Engel cone example}\label{sengel}

In this section we revisit the example in \cite{fsc} of a set with a
constant normal which is not a vertical halfspace, and we show why
the improvement procedure does not work, at least at some points, in
this case.

\subsection{The Engel group}
Let us recall the definition of Engel Lie algebra and group.

Let $\E$ be the Carnot group whose Lie algebra is $\g = V_1\oplus
V_2 \oplus V_3$ with $V_1 = \span \{X_1 , X_2 \}$, $V_2 = \{\R X_3
\}$ and $V_3 =\{\R X_4 \}$, the only non zero commutation relations
being
\begin{equation}\label{engel5}
[X_1 , X_2 ] = -X_3 ,\qquad [X_1 , X_3 ] = -X_4 .
\end{equation}
An explicit representation of the vector fields in $\R^4$ is:
\begin{eqnarray*}
&&X_1 = \partial_1,\\
&&X_2 = \partial_2 -x_1 \partial_3 + \frac{x^2_1}{2} \partial_4 ,\\
&&X_3 = \partial_3 - x_1 \partial_4 ,\\
&&X_4 = \partial_4 .
\end{eqnarray*}

Clearly $\E$ is a Carnot group with step $s=3$, topological
dimension $n=4$, homogeneous dimension $Q=2\cdot 1+1\cdot 2+1\cdot
3=7$, and dimension of the horizontal layer $m=2$. From now on, we
shall use the coordinates above to denote the elements of the
group.

\subsection{A cone in the Engel group}

For any $\alpha>0$, let $P=P_\alpha :\R^4\to\R$ be the polynomial
$$P(x) = \alpha x_2^3 + 2x_4 ,$$
whose gradient is
$$
\nabla P(x)=\left(0, 3 \alpha x_2^2, 0 ,2\right).
$$
In particular all level sets $\{P=c\}$ of $P$ are obviously graphs
of smooth functions depending on $(x_1,x_2,x_3)$. The derivative of
$P$ is particularly simple along the vector fields of the horizontal
layer: indeed, we have
\begin{eqnarray*}
X_1P(x)&=&\partial_1 (\alpha x_2^3 + 2x_4)=0
\end{eqnarray*}
and
\begin{eqnarray*}
X_2P(x)&=&[ \partial_2 -x_1 \partial_3 + \frac{x^2_1}{2} \partial_4] (\alpha x_2^3 + 2x_4)\\
&=&  3 \alpha x_2^2+x_1^2\geq0.
\end{eqnarray*}
Hence
\begin{equation}\label{engel1}
X_1P(x)=0,\quad X_2P(x)=x_1^2+3\alpha x_2^2\qquad\forall
x\in\R^4.
\end{equation}
We define
$$ C :=\{x \in\R^4 :\ P (x) \leq 0\},$$
whose boundary $\partial C$ is the set $\{P=0\}$. Notice that, due
to the (intrinsic) homogeneity of degree 3 of the polynomial, the
set $C$ is a cone, i.e. $\delta_rC=C$ for all $r>0$.

We shall denote by $\nu^{eu}_C(x)=\nabla P(x)/|\nabla P(x)|$ the unit (Euclidean)
outer normal to $C$. We also have the expansion
\begin{eqnarray*}\label{engel3}
|\nabla P|(x) &=&\sqrt{4+9\alpha^2x_2^4}\\
&=& 2+\frac{9}{2}\alpha^2x_2^4+O(d^4(x,0)).
\end{eqnarray*}
Thanks to Subsection \ref{rgauss} the set $C$ has locally finite
perimeter, and more precisely we have the formula (\ref{Euclid}) (throughout this section ${\mathscr
H}^k$ is the Hausdorff measure induced by the Euclidean distance)
\begin{equation}\label{engel2}
Z\uno_C=-\frac {Z P}{|\nabla P|}{\mathscr
H}^3\llcorner_{\partial C} \qquad\forall Z\in\g.
\end{equation}
In particular \eqref{engel1} and \eqref{engel2} give
$$
D\uno_C=(X_1\uno_C,X_2\uno_C)=(0,1)X_2\uno_C=-\frac{x_1^2+3\alpha x_2^2}{|\nabla
P(x)|}(0,1){\mathscr H}^3\llcorner_{\partial C}.
$$
It follows that
\begin{equation}\label{engel4}
|D\uno_C|=\frac{x_1^2+3\alpha x_2^2}{|\nabla P(x)|}{\mathscr
H}^3\llcorner_{\partial C}
\end{equation}
and that the horizontal normal, that is the vector field $\nu_C=(0,1)$, is constant,
so that
\emph{all} points of ${\rm supp\,}|D\uno_C|$ belong to ${\mathcal
F}C$.

Since we proved in Lemma~\ref{provafis} that non-horizontal regular
directions $Z$ for $E$ give rise, after blow-up, to invariant
directions, at least at points $\bar x$ where $|Z\uno_E|(B_r(\bar
x))/r^{Q-2}$ is infinitesimal as $r\downarrow 0$, and since the cone
is self-similar under blow-up at $\bar x=0$, it must happen that
$|Z\uno_C|(B_r(0))/r^{Q-2}$ is not infinitesimal as $r\downarrow 0$
for any non-horizontal regular directions $Z$ (actually, for the
cone $C$, \emph{all} directions are regular). Let us show explicitly
this fact for $Z:=\Ad_{\exp(X_1)}(X_2)$: taking into account the
commutator relations \eqref{engel5} and
\begin{eqnarray*}
Z:=\Ad_{\exp(X_1)}X_2&=&X_2+[X_1,X_2] +  \frac{1}{2}[X_1,[X_1,X_2]]\\
&=&X_2-X_3+ \frac{1}{2}X_4\\
&=& \partial_2 -x_1 \partial_3 + \frac{x^2_1}{2} \partial_4- \partial_3 + x_1 \partial_4 +\partial_4\\
&=&\partial_2 - (1+ x_1 )\partial_3+\frac{(1+x_1)^2}{2}\partial_4.
\end{eqnarray*}
We can now compute the derivative along the vector field $Z$:
\begin{eqnarray*}
ZP(x)&=&\left[\partial_2 - (1+ x_1 )\partial_3+\frac{(1+x_1)^2}{2}\partial_4\right]( \alpha x_2^3 + 2x_4 )\\
&=& 3 \alpha x_2^2  + (1+x_1)^2 \\
&=&1+O(d(x,0)).
\end{eqnarray*}

Intuitively, the quotient $|Z\uno_C|(B_r(0))/|D\uno_C|(B_r(0))$
tends to $+\infty$ as $r\downarrow 0$ because of the relations
\eqref{engel2} and \eqref{engel4}, and the fact that $ZP(0)\neq 0$
(notice that the factor $|\nabla P|$ is close to $2$ near to the
origin). Let us make a more precise analysis: according to the
ball-box theorem, balls $B_r(0)$ are comparable to the boxes
$$
Q_r:=[-r,r]^2\times [-r^2,r^2]\times [-r^3,r^3],
$$
so we will compute the density on these boxes, rather than on balls.
We shall assume, for the sake of simplicity, that $\alpha\in (0,2]$.
The homogeneity of $C$ and the fact that $0\in{\mathscr F}C$ give
$|D\uno_C|(Q_r)=c r^6$ for some positive constant $c$. The function
$$
x_4=-\alpha x_2^3:=g(x_1, x_2,x_3),
$$
whose graph is $\partial C$, has  absolute value strictly less than $r^3$, thus
$Q_r\cap\partial C$ is the graph of $g$
on the ``basis'' $[-r,r]^2\times [-r^2,r^2]$ of the box $Q_r$.  Moreover, since $g$ has
zero gradient at the origin,
\begin{eqnarray*}
{\mathscr H}^3(Q_r\cap\partial C)&=&\int_{[-r,r]^2\times [-r^2,r^2]}\sqrt{1+|\nabla g|^2} d\L^3\\
&\sim&\int_{[-r,r]^2\times [-r^2,r^2]} 1 d\L^3\\
&=& \L^3([-r,r]^2\times
[-r^2,r^2])=8 r^4.
\end{eqnarray*}
From \eqref{engel2} we obtain $|Z\uno_C|(Q_r)= 4r^4+o(r^4)$ and we
conclude that $|Z\uno_C|(B_r(0))/|D\uno_C|(B_r(0))\sim r^{-2}$.

\subsection{A counterexample to asymptotic stability of halfspaces}

Let us consider the set
\begin{equation}\label{toobad}
E:=\left\{(x_1,x_2,x_3,x_4)\in \E:\ x_2+\arctan (x_4)>0\right\}. 
\end{equation}
Since $X_1\uno_E=0$ and $X_2\uno_E\geq 0$ this set has a constant horizontal
normal, and clearly it is not an halfspace. On the other hand, it is not
difficult to check that the inclusions
$$
\left\{x_2>\frac{\pi}{2}\right\}\subseteq E\subseteq
\left\{x_2>-\frac{\pi}{2}\right\}
$$
imply that $E$ is asymptotic at infinity to the halfpspace $\{x_2>0\}$.

\subsection{Other constant normal sets in the Engel group}
We present here another family of sets that have constant horizontal
normal. This time we have a dependence on two parameters
$a,\,b\in\R$. Let $P_{a,b}:\R^4\to\R$ be the polynomial
$$P_{a,b}(x) = 2ax_4-bx_3+x_2.$$
Since $\partial_2P_{a,b}\neq0$, all level sets $\{P_{a,b}=c\}$ of $P_{a,b}$ are obviously  smooth manifolds.

Note that when both $a$ and $b$ are zero, the sub-level sets are vertical halfspaces.
In general, the derivatives along the vector fields of the horizontal
layer are
\begin{equation}\label{engelother}
X_1P(x)=0,\quad X_2P(x)=ax_1^2+b x_1+1\qquad\forall
x\in\R^4.
\end{equation}
So, if $(a,b)$ is close to $(1,0)$ then $ax_1^2+b x_1+1$ is a
perturbation of $x_1^2+1$ that is strictly greater than $0$. Thus
$X_2P(x)>0$ for any $(a,b)$ in a neighborhood of $(1,0)$. In other
words  the sub-level sets have constant horizontal normal. However,
these sets are not cones, except when they are vertical halfspaces.

\end{document}